\documentclass[a4paper,11pt]{amsart}
\usepackage{amsmath,amssymb,amsthm}
\usepackage{graphicx}
\usepackage{geometry}
\newtheorem{thm}{Theorem}[section]
\newtheorem{lemma}[thm]{Lemma}
\newtheorem{theorem}[thm]{Theorem}
\newtheorem{proposition}[thm]{Proposition}

\theoremstyle{definition}
\newtheorem{definition}[thm]{Definition}
\newtheorem{example}[thm]{Example}
\newtheorem{remark}[thm]{Remark}

\newtheorem{notation}[thm]{Notation}
\newtheorem{assumption}[thm]{Assumption}
\newtheorem{supplementation}[thm]{Supplementation}
\newtheorem{setting8}{Setting}

\newcommand{\RR}{\mathbb R}
\newcommand{\RRR}{\mathbb {R}_{\geq 0}}
\newcommand{\R}{\mathcal R}
\newcommand{\NN}{\mathbb N}

\newcommand{\HH}{\mathbb H}
\renewcommand{\H}{\mathcal H}
\newcommand{\ZZ}{\mathbb Z}
\newcommand{\F}{\mathcal F}
\newcommand{\Fla}{\mathcal{F}_{\lambda}}

\newcommand{\gL}{\Lambda}
\newcommand{\gG}{\Gamma}
\newcommand{\gS}{\Sigma}
\newcommand{\id}{\mathrm{id}}
\newcommand{\SCyl}{\mathrm{SCyl}}
\newcommand{\Cyl}{\mathrm{Cyl}}

\newcommand{\ti}[1]{\textit{#1}}

\renewcommand{\:}{\colon}
\renewcommand{\tilde}{\widetilde}
\renewcommand{\hat}{\widehat}
\newcommand{\SC}{\mathrm{SC}}
\newcommand{\Cay}{\mathrm{Cay}}
\newcommand{\Stab}{\mathrm{Stab}}

\newcommand{\GC}{\mathrm{GC}}

\newcommand{\ol}[1]{\overline{#1}}

\newcommand{\Comp}{\mathrm{Comp}}

\begin{document}
\title{Currents on cusped hyperbolic surfaces and denseness property}

\author[D. Sasaki]{Dounnu Sasaki}
\address{Department of Mathematics, Faculty of Science, Gakushuin University, Mejiro 1-5-1, Toshima-ku, Tokyo 171-8588, Japan}
\email{dounnu-daigaku@moegi.waseda.jp}

\begin{abstract}
The space $\mathrm{GC}(\Sigma)$ of geodesic currents on a hyperbolic surface $\Sigma$ can be considered as a completion of the set of weighted closed geodesics on $\Sigma$ when $\Sigma$ is compact, since the set of rational geodesic currents on $\Sigma$, which correspond to weighted closed geodesics, is a dense subset of $\mathrm{GC}(\Sigma )$.
We prove that even when $\Sigma$ is a cusped hyperbolic surface with finite area, $\mathrm{GC}(\Sigma )$ has the denseness property of rational geodesic currents, which correspond not only to weighted closed geodesics on $\Sigma$ but also to weighted geodesics connecting two cusps.
In addition, we present an example in which a sequence of weighted closed geodesics converges to a geodesic connecting two cusps, which is an obstruction for the intersection number to extend continuously to $\mathrm{GC}(\Sigma )$.
To construct the example, we use the notion of subset currents. Finally, we prove that the space of subset currents on a cusped hyperbolic surface has the denseness property of rational subset currents.
\end{abstract}

\maketitle

\section{Introduction}

Let $\gS$ be a hyperbolic surface with finite area (possibly with geodesic boundary). Geodesic currents on $\gS$, which were introduced by Bonahon in \cite{Bon86} as a generalization of measured geodesic laminations on $\gS$, have been successfully studied when $\gS$ is closed or compact.
They have been employed in the study of the Teichm\"uller space, mapping class groups, Kleinian groups, counting curves problems, and so on (see \cite{EU18} for a recent survey).

A geodesic current on $\gS$ is a (positive) $\pi_1(\gS)$-invariant Radon measure on the space
\[ \partial_2^\infty \tilde{\gS}:=\{ S\subset \partial^\infty\tilde{\gS}\mid \#S =2\} \]
for the boundary at infinity, $\partial^\infty \tilde{\gS}$, of the universal cover $\tilde{\gS}$ of $\gS$. Note that an element of $\partial_2^\infty \tilde{\gS}$ corresponds to an unoriented geodesic line on $\tilde{\gS}$.
We endow the space $\GC(\gS)$ of geodesic currents on $\gS$ with the weak-$\ast$ topology.

For each closed geodesic $\gamma$ on $\gS$, we can define a counting geodesic current $\eta_\gamma$ by considering all the lifts of $\gamma$ to $\tilde{\gS}$.
We regard $c\eta_\gamma\ (c\in \RRR)$, which is called a rational geodesic current, as a weighted closed geodesic on $\gS$.
When $\gS$ is compact, it has been proven in \cite{Bon86} that the set of rational geodesic currents on $\gS$ is a dense subset of $\GC(\gS)$.
In this sense, we say that $\GC(\gS)$ has the denseness property (of rational geodesic currents).
For a general hyperbolic surface $\gS$ we say that $\GC(\gS)$ has the denseness property if the set of rational currents, which is a weighted ``discrete'' measure corresponding to the $G$-obits of some point of $\partial_2^\infty \tilde{\gS}$, is a dense subset of $\GC(\gS)$ (see Definition \ref{def:rational}).

However, when $\gS$ has cusps, it has not been proven that $\GC(\gS)$ has the denseness property.
We remark that if $\gS$ has cusps, then a geodesic $\ell$ connecting two cusps, which is the projection of a geodesic line connecting two parabolic fixed points of $\partial^\infty \tilde{\gS}$, also induces a counting geodesic current $\eta_\ell$ similarly by considering all the lifts of $\ell$.

In this paper,  we prove that the space $\GC(\gS)$ of geodesic currents on a cusped hyperbolic surface $\gS$ with finite area has the denseness property.
Our strategy for the proof is based on \cite{Bon91} and \cite{Sas19}.
For a given geodesic current $\mu \in \GC (\gS)$, we construct a $G$-invariant family of quasi-geodesics on $\tilde{\gS}$ that induces a sum of counting geodesic currents approximating $\mu$.
To use the method in the case of compact hyperbolic surfaces, we cut off cusps along horocyclic curves around the cusps. One of the aspects of the proof is that such horocyclic curves are chosen more closely to cusps as we approximate $\mu$ more precisely.
\smallskip

\noindent \textbf{Other results.}
In Section \ref{sec:convergence}, we present an example in which a sequence of weighted closed geodesics converges to any given geodesic connecting two cusps in the space $\GC(\gS)$ of geodesic currents on a cusped hyperbolic surface $\gS$.
As a result, we can see that the set of weighted closed geodesics is a dense subset of $\GC(\gS)$.
Moreover, we construct an example in which a sequence of weighted geodesics connecting two cusps converges to any given closed geodesic in $\GC(\gS)$, which implies that the set of weighted geodesics connecting two cusps is a dense subset of $\GC(\gS)$.

In Section \ref{sec:extension of intersection number}, for a cusped hyperbolic surface $\gS$, we present a concrete example to prove that the intersection number cannot extend continuously to $\GC(\gS)$. To construct the example, we use the sequence of weighted closed geodesics converging to a geodesic connecting two cusps in Section \ref{sec:convergence}.
Note that according to \cite[Theorem 2.4]{BIPP19}, if we restrict $\GC(\gS)$ to the subset consisting of ``compact supported'' geodesic currents, then the intersection number can be extended continuously.
We present another sketch of the proof by using the method in \cite{Sas19}, which was used to prove the continuity of the extension of the intersection number to the space of subset currents.

In Section \ref{sec:dense subset currents}, we prove that the space $\SC(\gS)$ of subset currents on a cusped hyperbolic surface $\gS$ also has the denseness property. The notion of subset currents was introduced by Kapovich and Nagnibeda as a natural generalization of geodesic currents, and the study of subset currents began with the case of subset currents on free groups in \cite{KN13} and \cite{Sas15}.
See \cite{Sas19} for the study of subset currents on a compact hyperbolic surface, where the denseness property of rational subset currents has been proven.

A subset current on a hyperbolic surface $\gS$ is a $\pi_1(\gS)$-invariant Radon measure on the hyperspace
\[ \H( \partial^\infty \tilde{\gS}):=\{ S\subset \partial^\infty\tilde{\gS}\mid S\:\text{closed}, \#S\geq 2\}.\]
For a finitely generated subgroup $H$ of $\pi_1(\gS)$ whose limit set $\gL(H)$ contains at least two points, which belongs to $\H( \partial^\infty \tilde{\gS})$, we can define a counting subset current $\eta_H\in \SC(\gS)$ by considering the $\pi_1(\gS)$-orbit of $\gL(H)$. Then the set of weighted counting subset currents is proven to be a dense subset of $\SC (\gS)$.
We remark that the property of $\SC(\gS)$ is quite different from that in the case of a compact hyperbolic surface because $\SC(\gS)$ includes more types of rational subset currents (see Theorem \ref{thm:characterize rational currents}). We present some interesting examples in Section \ref{sec:convergence}, one of which is useful for the proof of the denseness property of geodesic currents on $\gS$.
\medskip

\noindent \textbf{Acknowledgements.}
I would like to express my heartfelt gratitude to Prof. Katsu\-hiko Matsu\-zaki, who provided carefully considered feedback and valuable comments.
I also would like to thank the referee for her/his careful reading of the manuscript and valuable comments.
The author is partially supported by JSPS KAKENHI Grant Number JP19K14539 and Grant-in-Aid for JSPS Fellows 21J01271.

\section{Preliminary}

In this section, we present the definition of geodesic currents and subset currents in a unified manner, and we summarize the results of this paper.

Let $G$ be a group acting continuously on a topological space $X$. The topological spaces that we deal with in this paper are always locally compact, separable, and completely metrizable.

We denote by $M_G(X)$ the space of $G$-invariant (positive) locally finite Borel measures on $X$, and we endow $M_G(X)$ with the weak-$\ast$ topology, i.e., a sequence $\{ \mu_n\}$ of $M_G(X)$ converges to $\mu\in M_G(X)$ if and only if
\[ \int f d\mu_n\rightarrow \int f d\mu \quad (n\rightarrow \infty )\]
for every continuous function $f\: X\rightarrow \RR$ with compact support.
Recall that a Borel measure $\mu$ on $X$ is $G$-invariant if for any $g\in G$, the push-forward measure $g_\ast (\mu)$ of $\mu$ by $g$ is equal to $\mu$. Recall that $g_\ast(\mu )(E)=\mu (g^{-1}(E))$ for any Borel subset $E$ of $X$.
A Borel measure $\mu$ on $X$ is said to be locally finite if $\mu(K)$ is finite for any compact subset $K$ of $X$.
We note that a locally finite Borel measure on $X$, which is a locally compact Polish space, is inner regular and outer regular (see \cite[2.18 Theorem]{Rud86}), and hence satisfies the condition of a Radon measure.

\begin{definition}[Geodesic currents and subset currents on hyperbolic groups]
Let $G$ be an infinite (Gromov) hyperbolic group and let $\partial G$ be the (Gromov) boundary of $G$.
Set
\[ \partial _2G:=\{ S\subset \partial G\mid \#S =2\}\]
and
\[ \H(\partial G):=\{ S\subset \partial G\mid S\:\text{closed and }\# S\geq 2\}.\]
We endow $\H(\partial G)$ with the Vietoris topology, which is generated by the set of the forms
\[ \{ S\in \H( \partial G)\mid S\subset U\} \ \text{and}\ \{ S\in \H(\partial G) \mid S\cap U\not=\emptyset\}\]
for an open subset of $\partial G$.
This topology coincides with the topology induced by the Hausdorff distance on $\partial G$ with respect to some metric on $\partial G$ compatible with the topology.
We endow $\partial _2 G$ with the subspace topology of $\H (\partial G)$.
Note that the action of $G$ on $\partial G$ naturally induces the action of $G$ on $\partial _2G$ and $\H (\partial G)$.

We refer to $\GC (G):=M_G(\partial_2G)$ as the space of geodesic currents on $G$ and its elements as geodesic currents.
We refer to $\SC (G):=M_G(\H (\partial G))$ as the space of subset currents on $G$ and its elements as subset currents.
\end{definition}

\begin{definition}[Geodesic currents and subset currents on hyperbolic surfaces]
Let $\gS$ be a hyperbolic surface possibly with (closed) geodesic boundaries.
Hyperbolic surfaces that we deal with in this paper are always complete, oriented, and connected. In addition, we usually assume that a hyperbolic surface has finite area.
We consider the universal cover $\tilde{\gS}$ of $\gS$ as a convex subspace of the hyperbolic plane $\HH$.
Then the boundary at infinity of $\tilde{\gS}$, denoted by $\partial^{\infty}\tilde{\gS}$, is the limit set of $\tilde{\gS}$ in $\HH$, which is the set of accumulation points of $\tilde{\gS}$ in the boundary $\partial \HH$.
Note that the fundamental group $\pi_1(\gS)$ of $\gS$ acts on $\tilde{\gS}$ and on $\partial^\infty \tilde{\gS}$.
When $\pi_1(\gS)$ has a parabolic element $p$ as an isometry of $\HH$, the projection of the neighborhood of the fixed point $p^\infty$ of $p$ to $\gS$ is called a cusp neighborhood, and we call $\gS$ a \ti{cusped hyperbolic surface}.

For $\partial^\infty \tilde{\gS}$, we also use the notation $\partial^\infty_2 \tilde{\gS}$ and $\H (\partial ^\infty(\gS))$ that we have introduced above.
We refer to $\GC (\gS):=M_G(\partial^\infty_2 \tilde{\gS})$ as the space of geodesic currents on $\gS$ and its elements as geodesic currents.
We refer to $\SC (\gS):=M_G(\H (\partial^\infty \tilde{\gS}))$ as the space of subset currents on $\gS$ and its elements as subset currents.
\end{definition}

\begin{remark}[Motivation for this paper]\label{rem:motivation of paper}
If the hyperbolic surface $\gS$ is compact, then $G=\pi_1(\gS)$ is a hyperbolic group and there exists a natural $G$-equivariant homeomorphism $\phi_G: \partial G\rightarrow \partial^\infty \tilde{\gS}$, which means that the action of $G$ on $\partial G$ essentially equals the action of $G$ on $\partial^\infty \tilde{\gS}$.
Then we can see that $\phi_G$ induces the isomorphism from $\GC(G)$ to $\GC(\gS)$ and from $\SC(G)$ to $\SC (\gS)$.

However, when $\gS$ has some cusps, $G=\pi_1(\gS)$ is a free group of finite rank, which means that $G$ is a hyperbolic group; however, the property of $\GC(G)$ (or $\SC(G)$) is \ti{quite different} from those of $\GC (\gS)$ (or $\SC (\gS)$).
The main purpose of this paper is to investigate the spaces $\GC(\gS)$ and $\SC(\gS)$ in this case.

We remark that even when $\gS$ has some cusps, there exists a natural $G$-equivariant continuous map $\phi$ from $\partial G$ to $\partial^\infty \tilde{\gS}$, which is referred to as the Cannon--Thurston  map. However, $\phi$ is surjective but not injective. We consider this Cannon--Thurston map and its application in Section \ref{sec:Cannon-Thurston}.
\end{remark}

\begin{definition}
Let $G$ be a group acting continuously on a topological space $X$.
For $x\in X$, we define a $G$-invariant Borel measure $\eta_x
$ on $X$ as
\[ \eta_x=\sum_{g\Stab (x)\in G/\Stab (x)} \delta_{gx},\]
where $\Stab(x)=\{ g\in G\mid g(x)=x\}$ and $\delta_{gx}$ is the Dirac measure at $gx$ on $X$.
Note that $\eta_x(E)$ equals the number of $G$-orbits of $x$ for a Borel subset $E\subset X$.
\end{definition}

We remark that in the context of geodesic currents and subset currents, we need to see whether $\eta_S$ is locally finite or not for $S\in \partial^\infty_2\tilde{\gS}$ or $\H (\partial^\infty \tilde{\gS})$.
When the hyperbolic surface $\gS$ is compact, $\SC(\gS)$ is isomorphic to $\SC (\pi_1(\gS))$ and the following proposition solves this problem.

\begin{theorem}[See {\cite[Theorem 2.8]{Sas19}}]
Let $G$ be an infinite hyperbolic group.
Let $S\in \H (\partial G)$. The $G$-invariant measure $\eta_S$ is locally finite if and only if $H:=\Stab (S)$ is a quasi-convex subgroup of $G$ and $S$ coincides with the limit set $\gL_G(H)$ of $H$. In particular, if a subset current $\mu \in \SC (G)$ has
an atom $S$, then $\Stab(S)$ is a quasi-convex subgroup of $G$ and $S=\gL_G(\Stab (S))$.
\end{theorem}

Note that if $S\in \partial_2G$, then the equality $S=\gL_G(H)$ for a subgroup $H$ of $G$ implies that $H$ is a subgroup $\langle h\rangle$ generated by $h\in G$.
When  $G$ is a free group of finite rank, a subgroup $H$ of $G$ is quasi-convex if and only if $H$ is finitely generated.

We generalize the above theorem to the case of subset currents on cusped hyperbolic surfaces and prove the following theorem (see Section \ref{sec:rational currents} for further details).
Recall that the limit set $\gL(H)$ of a subgroup $H$ of $\pi_1(\gS)$ (in $\tilde{\gS}$) is the limit set of the orbit $H(x)$ for some $x\in \tilde{\gS}$.

\begin{theorem}\label{thm:characterize rational currents}
Let $\gS$ be a cusped hyperbolic surface with finite area and let $G$ be the fundamental group of $\gS$.
Let $S\in \H (\partial^\infty \tilde{\gS})$.
The $G$-invariant measure $\eta_S$ is locally finite if and only if $H:=\Stab(S)$ is a finitely generated subgroup of $G$, and there exists a finite set $P$ of parabolic fixed points of $\partial^\infty \tilde{\gS}$ (possibly empty) such that $S=\gL(H)\sqcup H(P)$.
Note that $H$ can be a trivial subgroup $\{ \id \}$; then, $P$ contains at least two points.
\end{theorem}

\begin{remark}
Let $\gS$ be a cusped hyperbolic surface with finite area and let $G$ be the fundamental group of $\gS$.
For a non-trivial finitely generated subgroup $H$ of $G$ with $\#\gL(H)\geq 2$, we can consider the $G$-invariant measure $\eta_{\gL(H)}$. Then we need to note that $\Stab(\gL(H))$ is not necessarily equal to $H$. In general, $H$ is a finite-index subgroup of $\Stab(\gL (H))$, which implies that $\Stab(\gL(H))$ is also finitely generated. Therefore, $\eta_{\gL(H)}$ is locally finite.
We define
\[ \eta_H:=\sum_{gH\in G/H} \delta_{g\gL(H)}.\]
Then we see that $\eta_H=k\eta_{\gL(H)}$ if $H$ is a $k$-index subgroup of $\Stab (H)$.

When $H=\langle g\rangle$ for a hyperbolic element $g\in G$, we write $\eta_g$ in place of $\eta_{\langle g\rangle}$.

When $\#\gL(H)=1$, i.e., $H=\langle g\rangle$ for a parabolic element $g\in G$, we consider $\eta_H$ as the zero measure on $\H (\partial^\infty \tilde{\gS})$ for convenience.
In Theorem \ref{thm:characterize rational currents}, if $\#\gL(H)=1$, then $P$ contains at least one point.
\end{remark}

\begin{definition}\label{def:rational}
Let $G$ be a group acting continuously on a topological space $X$.
We say that $\mu \in M_G(X)$ is \ti{rational} if there exist $c\geq 0$ and $x\in X$ such that $\mu =c\eta_x$, and set
\[ M_G^r(X)=\{ \mu\in M_G(X)\mid \mu \: \text{rational}\}.\]
We say that $\mu \in M_G(X)$ is \ti{discrete} if there exist $c_1,\dots ,c_n\geq 0$ and $x_1,\dots ,x_n \in X$ such that
\[ \mu =c_1\eta_{x_1}+ \cdots +c_n\eta_{x_n},\]
and set
\[ M_G^d(X)=\{ \mu \in M_G(X)\mid \mu \:\text{discrete}\},\]
which is an $\RRR$-linear span of $M_G^r(X)$. Note that the discrete $G$-invariant measure $\mu \in M_G^d(X)$ is different from the usual discrete measure, which can be an infinite sum of rational measures.
\end{definition}

In the context of geodesic currents and subset currents, we will use the notation
\[\GC^r(G),\GC^d(G), \SC^r(G), \GC^r (\gS ),...\]
and so on to denote rational or discrete currents.
We summarize some important theorems related to the denseness property of such subsets.

\begin{theorem}[See \cite{Bon91}]\label{thm: hyp gp gc dense}
Let $G$ be an infinite hyperbolic group.
The set $\GC^r(G)$ of rational geodesic currents is a dense subset of $\GC(G)$.
\end{theorem}

\begin{theorem}[See \cite{KN13}]\label{thm:free gp dense}
Let $F$ be a free group of finite rank.
The set $\SC^r(F)$ of rational subset currents is a dense subset of $\SC (F)$.
\end{theorem}

\begin{theorem}[See \cite{Sas19}]
Let $G$ be a surface group, i.e., the fundamental group of a closed hyperbolic surface.
The set $\SC^r(G)$ of rational subset currents is a dense subset of $\SC (G)$.
\end{theorem}

We remark that if we consider the case in which a (discrete) group $G$ acts on $X$ properly discontinuously, the denseness of $M_G^r(X)$ in $M_G(X)$ might seem to be unnatural.
However, the action of a hyperbolic group $G$ on $\partial_2 G$ and on $\H (\partial G)$ is far from properly discontinuous, and we can prove that $\GC^r(G)$ is a dense subset of $\GC^d (G)$.

We also remark that for a general infinite hyperbolic group $G$, it is still open whether $\SC^r(G)$ is a dense subset of $\SC^d (G)$ or not.

From the above-mentioned results, it is natural to consider the question of whether the space of geodesic currents (or subset currents) on a ``cusped'' hyperbolic surface has such a denseness property.
The following two theorems are the main results of this paper.

\begin{theorem}\label{thm:GC denseness 1}
Let $\gS$ be a cusped hyperbolic surface with finite area.
The set $\GC^r(\gS)$ of rational geodesic currents is a dense subset of $\GC (\gS)$.
\end{theorem}

\begin{theorem}\label{thm:SC denseness 1}
Let $\gS$ be a cusped hyperbolic surface with finite area.
The set $\SC^r(\gS)$ of rational subset currents is a dense subset of $\SC (\gS)$.
\end{theorem}

From the results of Theorem \ref{thm:dense of subgp and current}, we can improve Theorems \ref{thm:GC denseness 1} and \ref{thm:SC denseness 1} as follows:

\begin{theorem}\label{thm:GC dense closed geodesic}
Let $\gS$ be a cusped hyperbolic surface with finite area.
The set
\[ \{ c\eta_g \mid c>0,\ g\in \pi_1(\gS) \: \text{hyperbolic element}\}\]
is a dense subset of $\GC (\gS)$. Note that $c\eta_g$ corresponds to a weighted closed geodesic on $\gS$.
\end{theorem}

\begin{theorem}
Let $\gS$ be a cusped hyperbolic surface with finite area.
The set
\[ \{ c\eta_H \mid c>0,\ H< \pi_1(\gS) \: \text{finitely generated subgroup}\}\]
is a dense subset of $\GC (\gS)$.
\end{theorem}

In addition, we can also obtain the following denseness theorem as a corollary of Proposition \ref{prop:geodesic line converges hyperbolic element}.

\begin{theorem}\label{thm:dense of geodesic line of para}
Let $\gS$ be a cusped hyperbolic surface with finite area.
The set
\[ \{ c\eta_{ \{ p,q \} } \mid c>0,\ p,q\in \partial^\infty \tilde{\gS}\: \text{disjoint parabolic fixed points} \}\]
is a dense subset of $\GC (\gS)$. Note that $c\eta_{\{ p,q\}}$ corresponds to a geodesic connecting two cusps.
\end{theorem}

\section{Rational currents on cusped hyperbolic surfaces}\label{sec:rational currents}

Let $\gS$ be a cusped hyperbolic surface with finite area and let $G$ be the fundamental group of $\gS$.
In this section, we present the proof of Theorem \ref{thm:characterize rational currents}.

\begin{assumption}\label{assump:gS has no boundary}
In general, $\gS$ can have some geodesic boundary; however, for simplicity, we assume $\gS$ has no boundary throughout the paper.
Then we regard the hyperbolic plane $\HH$ as the universal cover of $\gS$.
Actually, in most cases, the same argument works for a cusped hyperbolic surface with boundary by replacing $\HH$ with the universal cover $\tilde{\gS}$ of $\gS$.
Let $\pi$ be the canonical projection from $\HH$ to $\gS$.
\end{assumption}

\begin{definition}[Horocycle parameter]
For each cusp $p$ of $\gS$, we can take a parabolic fixed point $\xi \in \partial \HH$ corresponding to $p$.
A horocycle $h$ around $p$ is the projection of a horocycle centered at $\xi$ from $\HH$ to $\gS$.

Assume that $\gS$ has $k$ cusps and take disjoint horocycles $h_1,\dots ,h_k$ around each cusp.
Then we refer to a set $\{ h_1,\dots , h_k\}$ of horocycles as a \ti{horocycle parameter} of $\gS$.
We say that a horocycle parameter $\lambda=\{ h_1,\dots ,h_k\}$ is \ti{large} if each horocycle $h_i$ around a cusp $p_i$ is ``close'' to $p_i$ for $i=1,\dots ,k$, i.e., the radius of a horocycle centered at $\xi_i$ corresponding to $h_i$ is small.

For a horocycle parameter $\lambda=\{ h_1,\dots ,h_k\}$, we define $\gS_\lambda$ as a surface with non-geodesic boundaries obtained by cutting off cusps from $\gS$ along each $h_i$. We assume that $\gS_\lambda$ includes the horocycles $h_1,\dots ,h_k$, which implies that $\gS_\lambda$ is a compact subset of $\gS$.
\end{definition}

Fix some horocycle parameter $\lambda$.
We set $\HH_\lambda:=\pi^{-1}(\gS_\lambda) \subset \HH$.
Take a Dirichlet fundamental domain $\F$ corresponding to the action of $G$ on $\HH$.
Then we see that $\F_\lambda:= \F \cap \HH_\lambda$ is a compact fundamental domain corresponding to the action of $G$ on $\HH_\lambda$.

Recall that for $S\in \H (\partial \HH)$, the \ti{convex hull} $CH(S)$ of $S$ is the smallest convex subset of $\HH$ including all geodesic lines connecting two points of $S$.
For a bounded subset $K$ of $\HH$, we define a set $A(K)$ as
\[ A(K):=\{ S\in \H (\partial \HH )\mid CH(S)\cap K\not= \emptyset\}.\]
From \cite[Lemma 3.7 and 3.8]{Sas19}, $A(K)$ is a relatively compact subset of $\H (\partial \HH)$, and for any compact subset $E$ of $\H (\partial \HH)$, there exists a bounded subset $K$ of $\HH$ such that $A(K)$ includes $E$. Moreover, if $K$ is compact, so is $A(K)$.

\begin{lemma}\label{lem:locally finite A(K) finite}
Let $\mu$ be a $G$-invariant Borel measure on $\H (\partial \HH)$.
The measure $\mu $ is locally finite if and only if $\mu (A(\F_\lambda ))<\infty $ for some $\lambda$.
\end{lemma}
\begin{proof}
The ``only if'' part follows immediately since $\F_\lambda$ is a compact subset of $\HH$ and $A(\F_\lambda)$ is a compact subset of $\H (\partial \HH)$.

We prove the ``if'' part. Take any compact subset $E$ of $\H (\partial \HH)$. Then we can take a compact subset $K$ of $\HH$ such that $A(K)$ includes $E$.
For $K\cap \HH_\lambda$, there exist $g_1,\dots ,g_m\in G$ such that
\[ K\cap \HH_\lambda \subset g_1\F_\lambda \cup \cdots \cup g_m \F_\lambda.\]
Note that $\HH \setminus \HH_\lambda$ is a union of infinite open horodisks, and $K$ intersects at most finitely many open horodisks $H_1,\dots ,H_n$.
Hence, we have
\[ K\subset  g_1\F_\lambda \cup \cdots \cup g_m \F_\lambda \cup (K\cap H_1)\cup \cdots \cup (K\cap H_n),\]
which implies that
\[ A(K)\subset A(g_1\F_\lambda )\cup \cdots \cup A(g_m \F_\lambda )\cup A(K\cap H_1)\cup \cdots \cup A(K\cap H_n).\]
Since $A(g_i\F_\lambda)=g_iA(\F_\lambda)$ for $i=1,\dots ,m$, we see that
\[ \mu (A(g_i\F_\lambda))=\mu (g_iA(\F_\lambda ))=\mu (A(\F_\lambda ))<\infty.\]
Hence, it is sufficient to see that $\mu (A(K)\cap H_j)<\infty $ for $j=1,\dots ,n$.

Consider the upper-half plane model of $\HH$ and assume that $H_j$ is a horodisk centered at $\infty$, i.e.,
\[ H_j=\{ x+iy\in\HH \mid x,y\in \RR, t<y\} \]
for some $t>0$.
Since $K$ is compact, we can take $a,b,c\in \RR$ such that
\[ K\cap H_j\subset \{ x+iy\in \HH\mid a\leq x\leq b, t<y\leq c \}.\]
Moreover, there exist $a',b'\in \RR$ such that if a geodesic line $\ell $ on $\HH$ intersects $K\cap H_j$, then $\ell $ must intersect the segment
\[ \{ x+it\in \HH \mid a'\leq x\leq b'\}.\]
The point is that this segment can be covered by a finite union of $g \F_\lambda\ (g\in G)$.
Therefore, $\mu (A(K))$ is finite and so is $\mu (E)$, which implies that $\mu$ is locally finite.
\end{proof}

From the argument in the above proof, we see that for any geodesic line $\ell$ on $\HH$, there exists $g\in G$ such that $\ell$ intersects $g\F_\lambda$.
Hence,
\[ G(A(\F_\lambda ))=\H (\partial \HH ) \text{ and } G(A(\F_\lambda )\cap \partial _2\HH )=\partial _2\HH ,\]
which implies that the actions of $G$ on $\partial_2 \HH$ and on $\H (\partial \HH )$ are cocompact.
From \cite[Theorem 2.23]{Sas19}, we can obtain the following proposition.

\begin{proposition}\label{prop:topology of GC and SC}
Let $\gS$ be a cusped hyperbolic surface with finite area.
The space $\GC (\gS)=M_G(\partial^\infty_2 \tilde{\gS})$ of geodesic currents on $\gS$ and the space $\SC (\gS)=M_G(\H (\partial^\infty \tilde{\gS}))$ of subset currents on $\gS$ are locally compact, separable, and completely metrizable spaces.
\end{proposition}

From the above proposition and Theorem \ref{thm:GC dense closed geodesic}, we can regard $\GC(\gS)$ as a ``completion'' of weighted closed geodesics on $\gS$.

We apply Lemma \ref{lem:locally finite A(K) finite} to the $G$-invariant Borel measure $\eta_S$ for some $S\in \H (\partial \HH)$.
Set $H=\Stab (S)$.
Then
\begin{align*}
\eta_S(A(\F_\lambda) )
&=\# \{ gH\in G/H\mid gS\in A(\F_\lambda )\} \\
&=\# \{ gH\in G/H\mid g CH(S)\cap \F_\lambda \not=\emptyset \}\\
&=\# \{ gH\in G/H\mid g (CH(S)\cap \HH_\lambda )\cap \F_\lambda \not=\emptyset \}\\
&=\# \{ gH\in G/H\mid (CH(S)\cap \HH_\lambda )\cap g^{-1}\F_\lambda \not=\emptyset \}.
\end{align*}
To count the number of cosets $gH$, we consider a fundamental domain corresponding to the action of $H$ on $CH(S)\cap \HH_\lambda$.
Then we can obtain the following lemma.

\begin{lemma}
Let $S\in \H (\partial \HH )$ and $H=\Stab(S)$. Let $\F_S$ be a Dirichlet fundamental domain corresponding to the action of $H$ on $CH(S)$.
Note that $H$ can be $\{ \id \}$; then, $\F_S=CH(S)$.
The measure $\eta_S$ is locally finite if and only if $\HH_\lambda \cap \F_S$ is compact for some $\lambda$.
\end{lemma}

\begin{proof}
The point of this proof is the local finiteness of a Dirichlet fundamental domain, i.e., any compact subset $K\subset \HH$ intersects only finitely many translates of a Dirichlet fundamental domain.

First, we consider the ``if'' part.
For $gH\in G/H$, assume that $(CH(S)\cap \HH_\lambda)\cap g^{-1}\F_\lambda\not=\emptyset$. Then there exists $h\in H$ such that
\[ h(\F_S\cap \HH_\lambda) \cap g^{-1}\F_\lambda\not=\emptyset .\]
This implies that $(\F_S\cap \HH_\lambda)\cap h^{-1}g^{-1}\F_\lambda\not=\emptyset$.
Hence, the number of such $h^{-1}g^{-1}\in G$ is at most finite since $\F_S\cap \HH_\lambda$ is compact. Therefore, $\eta_S(A(\F_\lambda))$ is finite, which implies that $\eta_S$ is locally finite.

Now, we consider the ``only if'' part. We prove the contraposition. Assume that $\HH_\lambda \cap \F_S$ is not compact. Then there exists an infinite distinct sequence $\{ g_i\} _{i\in \NN}$ of $G$ such that
\[ (\F_S\cap \HH_\lambda ) \cap g_i^{-1}\F_\lambda \not=\emptyset.\]
Note that if $g_iH=g_jH$ for $i\not=j$, then there exists $h\in H$ such that $g_ih=g_j$, and we have
\[ (\F_S\cap \HH_\lambda ) \cap g_i^{-1}\F_\lambda \not=\emptyset \text{ and } h(\F_S\cap \HH_\lambda ) \cap g_i^{-1}\F_\lambda \not=\emptyset.\]
Since $g_i\F_\lambda$ is compact, the number of such $h\in G$ is at most finite. Hence, $\# \{ g_iH \mid i\in \NN\} =\infty$, which implies that  $\eta_S(A(\F_\lambda))=\infty$.
\end{proof}

From the above lemma, we see that the compactness of $\HH_\lambda \cap \F_S$ is independent of the base point of the Dirichlet fundamental domain $\F_S$ and the horocycle parameter $\lambda$.

\begin{lemma}
Let $S\in \H (\partial \HH )$ and $H=\Stab(S)$. Let $\F_S$ be a Dirichlet fundamental domain corresponding to the action of $H$ on $CH(S)$.
The intersection $\F_S \cap \HH_\lambda$ is compact if and only if $H$ is finitely generated and there exists a finite set $P$ of parabolic fixed points of $\partial \HH$ (possibly empty) such that $S=\gL(H)\sqcup H(P)$.
Note that $H$ can be a trivial subgroup $\{ \id \}$; then, $P$ contains at least two points.
\end{lemma}

\begin{proof}
When $S$ consists of two points, the statement follows immediately. Hence, we assume that $\# S \geq 3$.

When $H$ is trivial, $S=P$ and $\F_S=CH(S)$; then the statement follows immediately.

When $H$ is generated by one parabolic element $\alpha$ of $G$, then $\gL(H)=\{ \alpha^\infty \}$ and $CH(\gL(H))$ is empty.
Then for an appropriate base point for the Dirichlet domain $\F_S$ we can take $P\subset S$ such that $S=\gL(H)\sqcup H(P)$ and $\F_S$ is a convex hull of $\gL(H)\sqcup P$. Note that $P$ contains a point $x$ such that $\alpha(x)\in P$. Hence, $\F_S \cap \HH_\lambda$ is compact if and only if $P$ is a (non-empty) finite set of parabolic fixed points.

We assume that $\#\gL(H)\geq 2$ hereafter.
The quotient space $CH(\gL(H))/H$, which is called the convex core of $H$, has finite area or is a circle if and only if $H$ is finitely generated.
Let $\F_H$ be a Dirichlet fundamental domain corresponding to the action of $H$ on $CH(\gL(H))$.

First, we prove the ``if'' part. Since $H$ is finitely generated, $\F_H$ has finite area. Hence, $\HH_\lambda \cap \F_H$ is compact. Moreover, the quotient space of the action of $H$ on each connected component of $CH(S)\setminus CH(\gL (H))$, which is called a \ti{crown}, also has finite area.
Therefore, $\HH_\lambda \cap \F_S$ is compact.

Next, we prove the ``only if'' part. We can assume that the base point of $\F_S$ coincides with the base point of $\F_H$.
Since $\HH_\lambda \cap \F_S$ is compact, so is $\HH_\lambda \cap \F_H$, which implies that $H$ is finitely generated.
Moreover, $\F_S$ and $\F_H$ are finite polygons whose vertices can be on $\partial \HH$.
If a vertex $v$ of $\F_S$ is on $\partial \HH$, then $v$ is a parabolic fixed point since $\HH_\lambda \cap \F_S$ is compact.
This implies that the set $P'$ of all vertices of $\F_S$ on $\partial \HH$ consists of finitely many parabolic fixed points.
Since $CH(S)=H(\F_S)$, we can see that
\[ S=\overline{\gL(H)\cup H(P')}.\]
Set $P=P'\setminus \gL(H)$.
Since $H$ acts on $\partial \HH \setminus \gL(H)$ properly discontinuously, we have
\[ S=\gL(H) \sqcup H(P),\]
as required.
\end{proof}

From the above lemmas, Theorem \ref{thm:characterize rational currents} follows.

\begin{remark}\label{rem:discrete subset current}
In the proof of the denseness property, we will construct $\nu \in \SC (\gS)$ denoted as
\[ \nu =\sum_{i\in I} \delta_{S_i}\]
for $S_i\in \H (\partial^\infty \tilde{\gS})\ (i\in I)$.
Since $\nu$ is $G$-invariant, $\eta_{S_i}$ is also a subset current on $\gS$ for every $i\in I$.
Moreover, there are finite $i_1,\dots i_k\in I$ such that
\[ \nu =\eta_{S_{i_1}}+\cdots +\eta_{S_{i_k}}\]
by Lemma \ref{lem:locally finite A(K) finite}, which implies that $\nu$ is a discrete subset current on $\gS$.
\end{remark}

\section{Cannon--Thurston maps and currents}\label{sec:Cannon-Thurston}

Let $\gS$ be a cusped hyperbolic surface with finite area and let $G$ be the fundamental group of $\gS$.
For simplicity, we assume that $\gS$ has no boundary; however, this assumption is not necessary.
Recall Remark \ref{rem:motivation of paper}.
We have the Cannon--Thurston map $\phi$ from $\partial G$ to $\partial^\infty \tilde{\gS}=\partial \HH$, which is a surjective continuous $G$-equivariant map sending $[ \{ x_n\} ]\in \partial G$ to the limit point of $x_n(y)\in \partial \HH$ for some point $y\in \HH$.
Then $\phi$ naturally induces a continuous map
\[ \H\phi\:\H (\partial G)\rightarrow \hat{\H} (\partial \HH )=\H (\partial \HH )\cup \{ \{ x \} \mid x\in \partial \HH \},\]
whose topology is the Vietoris topology, which coincides with the topology induced by a Hausdorff distance. Note that $S\in \H (\partial G)$ is a compact subset of $\partial G$ and so is $\phi (S)= \H \phi (S)$.

\begin{lemma}
For any compact subset $E$ of $\H (\partial \HH)$, the preimage $\H \phi^{-1}(E)$ is compact.
\end{lemma}
\begin{proof}
Take any compact subset $E$ of $\H (\partial \HH)$.
To obtain a contradiction, suppose that $\H\phi^{-1}(E)$ is not compact.
Note that $\hat{\H}(\partial G)$ is a compactification of $\H (\partial G)$ with respect to a Hausdorff distance.
Hence, we can take a sequence $S_n$ of $\phi^{-1}(E)$ converging to $\{ x\}$ for some $x\in \partial G$ in $\hat{\H} (\partial G)$.
Since $\H \phi$ is continuous, $\H\phi(S_n)$ converges to $\H\phi(\{ x\})$.
This implies that the compact set $E$ includes a sequence $\H\phi(S_n)$ converging to $\H\phi(\{ x \})\not\in E$, which is a contradiction.
\end{proof}

For $\mu \in \SC (G)$, by considering the push-forward by $\H \phi$, we have a $G$-invariant measure $\H\phi_\ast (\mu )$ on $\hat{\H}(\partial \HH)$ since $\H\phi$ is $G$-equivariant.
Then the restriction of the measure $\H\phi_\ast(\mu)$ to $\H (\partial \HH)$, denoted by $\phi_\SC(\mu )$, is locally finite from the above lemma.
Similarly, we can obtain a map $\phi_\GC$ from $\GC (G)$ to $\GC (\gS)$.
Note that the continuity of $\phi_\SC$ (and $\phi_\GC$) is not trivial since in the construction we restrict a measure on $\hat{\H} (\partial \HH)$ to $\H (\partial \HH)$.

\begin{lemma}
The maps $\phi_\SC\: \SC (G)\rightarrow \SC (\gS)$ and $\phi_\GC \: \GC (G)\rightarrow \GC (\gS)$ are $\RRR$-linear and continuous.
\end{lemma}
\begin{proof}
The $\RRR$-linearity follows immediately by the definition.
We prove that $\phi_\SC$ is continuous. The continuity of $\phi_\GC$ follows from the same proof.
Take a sequence $\mu_n\ (n\in \NN)$ of $\SC (G)$ converging to $\mu\in \SC (G)$.
It is sufficient to prove that $\phi_\SC (\mu_n)$ converges to $\phi_\SC (\mu)$.
Take any continuous function $f\: \H (\partial \HH)\rightarrow \RR$ with compact support.
Since we have \[\mathrm{supp}(f \circ \H\phi )\subset \H\phi^{-1}(\mathrm{supp}(f) ),\]
the support $\mathrm{supp}(f\circ \H\phi )$ is compact from the above lemma.
This implies that
\[ \int f d\phi_\SC(\mu_n)=\int f\circ \H\phi d \mu_n \underset{n\rightarrow \infty}{\longrightarrow} \int f\circ \H\phi d\mu=\int f d\phi_\SC (\mu).\]
Therefore, $\phi_\SC (\mu_n)$ converges to $\phi_\SC (\mu)$.
\end{proof}

By the definition of rational subset currents, we see that $\phi_\SC$ maps rational subset currents of $\SC (G)$ to rational subset currents of $\SC (\gS)$.
More concretely, for a non-trivial finitely generated subgroup $H$ of $G$, the limit set $\gL_G(H)$ of $H$ in $\partial G$ is mapped to the limit set $\gL(H)$ of $H$ in $\partial \HH$ by the map $\H\phi$. This implies that the subset current
\[ \eta_H^G=\sum_{gH\in G/H} \delta_{g\gL_G (H)} \in \SC (G) \]
is mapped to the subset current
\[ \eta_H=\sum_{gH\in G/H}\delta_{g\gL (H)} \in \SC (\gS)\]
by $\phi_\SC$. Note that if $H$ is a trivial subgroup, we define $\eta_H^G$ and $\eta_H$ to be the zero measure for convenience.

Similarly, for a non-trivial $h\in G$, the geodesic current
\[ \eta_h^G=\sum_{g\langle h\rangle G/\langle h\rangle} \delta_{g\gL_G (\langle h \rangle )}\in \GC (G)\]
is mapped to the geodesic current
\[ \eta_h=\sum_{g\langle h\rangle G/\langle h\rangle} \delta_{g\gL (\langle h \rangle )}\in \GC (\gS )\]
by $\phi_\GC$. Note that if $h$ is a parabolic element, then $\eta_h$ is the zero measure.

Recall that by Theorem \ref{thm:free gp dense} for a free group $F$ of finite rank,
\[ \SC^r(F) =\{ c\eta_H^G \mid c>0,\ H<G\:\text{finitely generated subgroup} \} \]
is a dense subset of $\SC^d(F)=\mathrm{Span}(\SC^r(F))$ (see Theorem \ref{thm: hyp gp gc dense} for the case of geodesic currents). Hence, by the continuity of $\phi_\SC$ and $\phi_\GC$, we can obtain the following lemma.

\begin{lemma}\label{lem:rational span dense}
The set
\[ \{ c\eta_H \in \SC (\gS )\mid c>0,\ H<G\:\text{finitely generated subgroup} \}.\]
is a dense subset of the $\RRR$-linear span
\[ \mathrm{Span}(\{ c\eta_H \in \SC (\gS )\mid c>0,\ H<G\:\text{finitely generated subgroup} \} ),\]
and the set
\[ \{ c\eta_h \in \GC (\gS )\mid c>0,\ h\in G\}\]
is a dense subset of the $\RRR$-linear span
\[ \mathrm{Span}(\{ c\eta_h \in \GC (\gS )\mid c>0,\ h\in G\}).\]
\end{lemma}

\begin{remark}
We remark that $\phi_\SC $ and $\phi_\GC $ are not surjective.
In fact, there exists no $S\in \H(\partial G)$ such that $\eta_S$ is rational, and $S$ is mapped to the set $\{\alpha^\infty ,\beta^\infty\}$ of two different parabolic fixed points for two parabolic elements $\alpha$ and $\beta$ of $G$.

If $\eta_S$ is rational, then for the stabilizer $H=\mathrm{Stab}(S)$, we have $S=\gL_G(H)$. Then $\phi(S)=\{ \alpha^\infty, \beta^\infty \}$ implies that $H$ includes $\alpha^k$ and $\beta^l$ for some $k,l\in \ZZ\setminus \{0\}$. Therefore, we have $S\supset \gL_G(\langle \alpha^k,\beta^l\rangle)$; then
\[ \{ \alpha ^\infty, \beta^\infty \} =\phi ( S ) = \gL(H) \supset \gL(\langle \alpha^k,\beta^l\rangle),\]
which is a contradiction.
\end{remark}

\section{Approximation of geodesic line by sequence of closed geodesics}\label{sec:convergence}

In this section, we present some interesting examples of convergence sequences of rational geodesic currents or subset currents on a cusped hyperbolic surface.
One of the examples is a sequence of closed geodesics converging to a weighted geodesic connecting two cusps, which is used for the proof of Theorem \ref{thm:GC denseness 1}.

Let $\gS$ be a cusped hyperbolic surface with finite area, and let $G$ be the fundamental group of $\gS$. Recall Assumption \ref{assump:gS has no boundary}.

\begin{proposition}\label{prop:converge seq}
Let $\alpha, \beta \in G$ be parabolic elements. Assume that the fixed point $\alpha^\infty \in \partial \HH$ of $\alpha$ is different from the fixed point $\beta^\infty \in \partial \HH$ of $\beta$.
Then the sequence $\{\eta_{\alpha^n\beta^n} \}$ converges to $2\eta_{\{\alpha^\infty, \beta^\infty\}}$ in $\GC (\gS)$ when $n$ tends to infinity. Note that $\alpha^n\beta^n$ is a hyperbolic element of $G$ for every $n\in \NN$.
\end{proposition}

We need the following lemma to prove the above proposition.

\begin{lemma}\label{lem:converge seq}
Let $\alpha, \beta \in G$ be parabolic elements with $\alpha^\infty \not= \beta^\infty$.
Then the sequence $\{\eta_{\langle \alpha^n,\beta^n\rangle} \}$ of subset currents converges to $\eta_{\{\alpha^\infty, \beta^\infty\}}$.
\end{lemma}
\begin{proof}
First of all, we note that the limit set $\gL (\langle \alpha^n ,\beta^n\rangle)$ of $\langle \alpha^n ,\beta^n\rangle$ converges to $\{ \alpha^\infty, \beta^\infty\}$ with respect to the Hausdorff distance on $\partial \HH$ by applying the Ping-Pong lemma to $\alpha^n$ and $\beta^n$ when $n$ tends to infinity.

Let $\F_n$ be a Dirichlet fundamental domain based at some point on the geodesic $[ \alpha^\infty ,\beta^\infty ]$ with respect to the action of $\langle \alpha^n ,\beta^n\rangle$ on the convex hull $CH (\langle \alpha^n ,\beta^n\rangle):=CH (\gL (\langle \alpha^n ,\beta^n\rangle))$.

Consider the compact subsurface $\gS_\lambda $ of $\gS$ with respect to some horocycle parameter $\lambda$. Set $\HH_\lambda=\pi^{-1}(\gS_\lambda)$.
Then we can see that for any horocycle parameter $\lambda$, the sequence $\F_n\cap \HH_\lambda$ of compact subsets converges to $[\alpha^\infty,\beta^\infty]\cap \HH_\lambda$ with respect to the Hausdorff distance on $\HH$.

Now, take any continuous function $f\: \H (\partial \HH )\rightarrow \RRR$ with compact support and take a compact subset $K$ of $\HH$ such that the support $\mathrm{supp}(f)$ of $f$ is included in $A(K)=\{ S\in \H (\partial \HH)\mid CH(S)\cap K\not=\emptyset \}$.
Take $r>0$ and assume that $\HH_\lambda \supset B(K,r)$.
Moreover, we assume that $n$ is large enough so that the Hausdorff distance between $\F_n\cap \HH_\lambda$ and $[\alpha^\infty,\beta^\infty]$ is smaller than $r$.

Set $G_0=\{ g\in G\mid g[\alpha^\infty,\beta^\infty]\cap B(K,r)\not=\emptyset \}$, which is a finite set
because $\eta_{\{\alpha^\infty,\beta^\infty\}}$ is a locally finite measure on $\partial_2\HH$ and $\eta_{\{\alpha^\infty,\beta^\infty\}}(A(B(K,r)))$ is finite.
We see that if $n$ is sufficiently large, then the map
\[ \phi \: G_0\rightarrow G/\langle \alpha^n,\beta^n\rangle; g\mapsto g\langle \alpha^n,\beta^n\rangle,\]
is injective. Actually, for $g_1,g_2 \in G_0$, if $g_1\langle \alpha^n,\beta^n\rangle =g_2\langle \alpha^n,\beta^n\rangle$, then
$g_2^{-1}g_1\in \langle \alpha^n,\beta^n\rangle $. Since $G_0$ is a finite set, so is $\{ g_2^{-1}g_1\mid g_1,g_2\in G_0\}$.
Note that $G=\pi_1(\gS)$ is a free group. There exists a largest positive integer $k$ such that $\alpha^{\pm k}$ or $\beta^{\pm k}$ appears in $\{ g_2^{-1}g_1\mid g_1,g_2\in G_0\}$ as a reduced word.
If $n$ is large than $k$, then $g_2^{-1}g_1=\id$, which implies that $g_1=g_2$.

By the definition, we have
\[ \int f d\eta_{\langle \alpha^n,\beta^n\rangle}=\sum_{g\langle \alpha^n,\beta^n\rangle\in G/\langle \alpha^n,\beta^n\rangle}f(g\gL (\langle \alpha^n,\beta^n\rangle ) ).\]
We want to prove that if $n$ is sufficiently large, then
\[ \sum_{g\langle \alpha^n,\beta^n\rangle\in G/\langle \alpha^n,\beta^n\rangle}f(g\gL (\langle \alpha^n,\beta^n\rangle ) )
= \sum_{g\in G_0}f(g\gL (\langle \alpha^n,\beta^n\rangle ) ).\]
Set
\[ J=\{ g\langle \alpha^n,\beta^n\rangle\in G/\langle \alpha^n,\beta^n\rangle\mid g\gL (\langle \alpha^n,\beta^n\rangle ) \in \mathrm{supp}(f)\}.\]
Then it is sufficient to prove that $\phi (G_0)\supset J$ for a large $n$.
Take any $g\gL (\langle \alpha^n,\beta^n\rangle )\in J$.
Then $gCH(\langle \alpha^n,\beta^n\rangle )\cap K\not=\emptyset$,
which implies that there exists $h\in \langle \alpha^n,\beta^n\rangle$ such that $gh\F_n\cap K\not=\emptyset$.
Since $\F_n\cap \HH_\lambda$ is included in the $r$-neighborhood of $[\alpha^\infty, \beta^\infty]\cap \HH_\lambda$, we have
\[ gh[\alpha^\infty, \beta^\infty ]\cap B(K,r)\not=\emptyset ,\]
which implies that $gh\in G_0$ and $\phi(gh)=g\langle \alpha^n,\beta^n\rangle $.

From the above, for a sufficiently large $n$, we have
\begin{align*}
\int f d\eta_{\langle \alpha^n,\beta^n\rangle}
&=\sum_{g\langle \alpha^n,\beta^n\rangle\in G/\langle \alpha^n,\beta^n\rangle}f(g\gL (\langle \alpha^n,\beta^n\rangle ) )\\
&=\sum_{g\in G_0}f(g\gL (\langle \alpha^n,\beta^n\rangle ) )\\
&\underset{n\rightarrow \infty}{\longrightarrow} \sum_{g\in G_0}f(g\{ \alpha^\infty, \beta^\infty \} )=\int f d\eta_{\{ \alpha^\infty,\beta^\infty\}}.
\end{align*}
This implies that $\eta_{\langle \alpha^n,\beta^n\rangle}$ converges to $\eta_{\{ \alpha^\infty,\beta^\infty\}}$.
\end{proof}

\begin{proof}[Proof of Proposition \ref{prop:converge seq}]
First, we note that $\langle \alpha^n\beta^n\rangle$ is a cyclic subgroup of  $\langle \alpha^n, \beta^n\rangle$.
Take any continuous function $f\: \partial_2\HH\rightarrow \RRR$ with compact support.
Then we can take $K\subset \HH, r>0$ and $G_0$ for $f$ as in the proof of Lemma \ref{lem:converge seq}.
We remark that from the argument in the above for a sufficiently large $n\in \NN$ and $g\langle \alpha^n,\beta^n \rangle\in (G/\langle \alpha^n,\beta^n\rangle) \setminus \phi(G_0)$, we see that
\[ gCH(\langle \alpha^n,\beta^n \rangle )\cap K=\emptyset, \]
which implies that
\[ gCH(\langle \alpha^n\beta^n\rangle )\cap K=\emptyset.\]

By using a bijection
\[ G/\langle \alpha^n\beta^n\rangle \rightarrow G/ \langle \alpha^n, \beta^n \rangle\times \langle \alpha^n,\beta^n\rangle/\langle \alpha^n\beta^n\rangle,\]
we have
\begin{align*}
\eta_{\langle \alpha^n\beta^n\rangle}
&=\sum_{g\langle \alpha^n\beta^n\rangle\in G/\langle \alpha^n\beta^n\rangle }\delta_{g\gL (\langle \alpha^n\beta^n\rangle )}\\
&=\sum_{g\langle \alpha^n,\beta^n\rangle\in G/\langle \alpha^n,\beta^n\rangle }\sum_{h\langle \alpha^n\beta^n\rangle\in \langle \alpha^n,\beta^n\rangle/\langle \alpha^n\beta^n\rangle }\delta_{gh\gL (\langle \alpha^n\beta^n\rangle )}.
\end{align*}
Hence, for a sufficiently large $n$, we see that
\[ \int fd\eta_{\langle \alpha^n\beta^n\rangle}
=\sum_{g\in  G_0}\sum_{h\langle \alpha^n\beta^n\rangle\in \langle \alpha^n,\beta^n\rangle /\langle \alpha^n\beta^n\rangle}f(gh\gL( \langle \alpha^n\beta^n\rangle) )\]

By considering the Ping-Pong of $\alpha^n$ and $\beta^n$ for a large $n$,
$gh\gL( \langle \alpha^n\beta^n\rangle )$ does not belong to the support of $f$ unless $h\langle \alpha^n\beta^n\rangle$ is equal to
\[ \langle \alpha^n\beta^n\rangle \ \mathrm{or}\ \alpha^{-n}\langle \alpha^n\beta^n\rangle =\beta^n \langle \alpha^n\beta^n\rangle .\]
We add a supplementary explanation to this claim after the proof.
Note that $\gL( \langle \alpha^n\beta^n\rangle )$ and $\alpha^{-n}\gL(\langle \alpha^n\beta^n\rangle)=\gL(\langle \beta^n\alpha^n\rangle)$ converge to $\{ \alpha^\infty,\beta^\infty \}$ when $n$ tends to infinity.
Therefore, for a sufficiently large $n$,
\begin{align*}
\int fd\eta_{\langle \alpha^n\beta^n\rangle}
&=\sum_{g\in  G_0}\Big( f (g\gL ( \langle \alpha^n\beta^n\rangle ))+f(g\gL (\langle \beta^n\alpha^n\rangle ) ) \Big)\\
&\underset{n\rightarrow \infty}{\rightarrow} \sum_{g\in G_0}2f(g\{ \alpha^\infty, \beta^\infty \} )=\int f d2\eta_{\{ \alpha^\infty,\beta^\infty\}}
\end{align*}
This implies that $\eta_{\langle \alpha^n\beta^n\rangle}$ converges to $2\eta_{\{ \alpha^\infty,\beta^\infty\}}$.
\end{proof}

We remark that $\langle \alpha^n\beta^n\rangle$ (or $\langle \alpha^{-n}\beta^n\rangle$) corresponds to the boundary of the convex core $CH(\langle \alpha^n,\beta^n\rangle )/\langle \alpha^n,\beta^n\rangle$. In this case, $\gL (\langle \alpha^n\beta^n\rangle )$ and $\gL (\langle \beta^n\alpha^n\rangle )$ correspond to the two boundary components of $CH(\langle \alpha^n,\beta^n\rangle )$, which converges to $[\alpha^\infty,\beta^\infty]$ with respect to the Hausdorff distance on $\ol{\HH}$ when $n$ tends to infinity.

From the proof of Lemma \ref{lem:converge seq}, we see that more generally $\eta_{\langle \alpha^m, \beta^n\rangle}$ converges to $\eta_{\{ \alpha^\infty,\beta^\infty\}}$ when $m$ and $n$ tend to infinity.
In addition, for finitely many parabolic elements $\alpha_1,\dots ,\alpha_k\in G (k\geq 2)$ whose fixed points $\alpha_1^\infty,\dots ,\alpha_k^\infty$ are pairwise distinct, we can see that
\[ \eta_{\langle \alpha_1^n , \dots ,\alpha_k^n \rangle} \underset{n\rightarrow \infty }{\longrightarrow} \eta_{\{ \alpha_1^\infty, \dots ,\alpha_k^\infty\}}. \]
More generally, we can prove the following proposition.

\begin{proposition}\label{prop:sc subgp approximate S}
Let $S\in \H (\partial \HH)$ satisfying the condition that $H=\mathrm{Stab}(S)$ is a finitely generated subgroup of $G$, and $S=\gL(H)\sqcup H(\{ \alpha_1^\infty, \dots ,\alpha_k^\infty\})$. Then we have
\[ \eta_{\langle H \sqcup \{ \alpha_1^n,\dots ,\alpha_k^n \} \rangle}\underset{n\rightarrow \infty}{\longrightarrow} \eta_S.\]
\end{proposition}
\begin{proof}
We present a sketch of the proof when $\gL(H)$ includes at least two points. Even in the other cases the following argument works.
Set $H_n=\langle H \sqcup \{ \alpha_1^n,\dots ,\alpha_k^n \} \rangle$.
By considering the Ping-Pong lemma, for a sufficiently large $n$, $H_n$ is the free product of $\langle \alpha_1^n\rangle,\cdots ,\langle \alpha_k^n\rangle$, and $H$.
Moreover, we see that the sequence of limit sets $\gL(H_n)$ converges to $S$ in $\H (\partial \HH)$.

Let $B$ be a boundary component of the convex hull $CH(H)=CH(\gL (H))$, which is a geodesic line connecting two points of $\gL(H)$.
Let $I$ be the open interval of $\partial \HH$ connecting the endpoints of $B$. We note that $I$ does not contain any points of $\gL(H)$ and assume that $I$ contains $\alpha_1^\infty,\dots, \alpha_k^\infty$.
Take some base point $b$ on $B$.
Take the Dirichlet fundamental domain $\F_n$ based at $b$ corresponding to the action of $H_n$ on $CH(H_n)$ and take the Dirichlet fundamental domain $\F$ based at $b$ corresponding to the action of $H$ on $CH(S)$.

From the choice of the base point $b$, we can see that for any horocycle parameter $\lambda$, the sequence $\F_n\cap \HH_\lambda$ of compact subsets converges to $\F\cap \HH_\lambda$ with respect to the Hausdorff distance on $\HH$.
(When $H$ is a trivial subgroup, the base point $b$ can be any point in $CH(H_n)$, and $\F=CH(S)$.)
Then by the same argument as in Lemma \ref{lem:converge seq}, we see that $\eta_{H_n}$ converges to $\eta_S$.
\end{proof}

From Lemma \ref{lem:rational span dense}, Proposition \ref{prop:converge seq} and Proposition \ref{prop:sc subgp approximate S},
we can obtain the following theorem.

\begin{theorem}\label{thm:dense of subgp and current}
Let $\gS$ be a cusped hyperbolic surface with finite area. The set
\[ \{ c\eta_H \in \SC (\gS )\mid c>0,\ H<\pi_1(\gS)\:\text{finitely generated subgroup} \}.\]
is a dense subset of
\[ \SC^d(\gS)=\{ \mu \in \SC (\gS )\mid \mu \:\text{discrete subset current} \} ,\]
and the set
\[ \{ c\eta_g \in \GC (\gS )\mid c>0,\ g\in \pi_1(\gS)\:\text{hyperbolic element}\}\]
is a dense subset of
\[ \GC^d(\gS)=\{ \mu \in \GC (\gS )\mid \mu \:\text{discrete geodesic current} \} .\]
\end{theorem}

Now, we prove the ``opposite'' of Proposition \ref{prop:converge seq}:

\begin{proposition}\label{prop:geodesic line converges hyperbolic element}
For every hyperbolic element $g\in G=\pi_1(\gS)$, there exists a sequence of pairs of parabolic fixed points $\{ p_n ,q_n\}$ of $\partial^\infty \tilde{\gS}$ such that $\frac{1}{2n} \eta_{\{p_n,q_n\}}$ converges to $\eta_g$ when $n$ tends to infinity.
\end{proposition}
\begin{proof}
Let $p,q$ be parabolic fixed points of $\partial \HH$ such that the geodesic line $[p,q]$ intersects with the axis $\mathrm{Ax}(g)=CH(\gL (\langle g \rangle ))=[g^\infty,g^{-\infty}]$ of $g$. We denote by $g^{-\infty}$ the repelling fixed point of $g$ and denote by $g^{\infty}$ the attracting fixed point of $g$.
Then we set
\[ p_n=g^{-n}p \text{ and } q_n=g^nq.\]
Note that $p_n$ converges to $g^{-\infty}$ and $q_n$ converges to $g^{\infty}$ when $n$ tends to infinity.
We prove that $\frac{1}{2n}\eta_{\{ p_n,q_n\}}$ converges to $\eta_g$. Our strategy is almost the same as in the proof of Proposition \ref{lem:converge seq}.

Take any continuous function $f\: \partial _2\HH \rightarrow \RRR$ with compact support and a compact subset $K\subset \HH$ such that $A(K)$ includes $\mathrm{supp}(f)$.
Let $x$ be the intersection point of $[p,q]$ and $\mathrm{Ax}(g)$.
For a positive integer $n$, set
\[ L_n:=[g^{-n}x, g^n x]\subset \mathrm{Ax}(g). \]
Define $P_n$ to be the geodesic ray from $g^{-n}x$ to $g^{-n}p$ and define $Q_n$ to be the geodesic ray from $g^n x$ to $g^n q$.
For a sufficiently large $n$, by combining $P_n$, $L_n$ and $Q_n$ we can obtain a quasi-geodesic line $\ell_n$ connecting $g^{-n}p$ to $g^nq$ which is included in the $r$-neighborhood of the geodesic line $[p_n,q_n]$ for some constant $r>0$ by the stability of geodesics (see Figure \ref{fig:converging-path}).
We can also assume that the $r$-neighborhood of $\ell_n$ includes $[p_n,q_n]$.

\begin{figure}[h]
\centering
\includegraphics{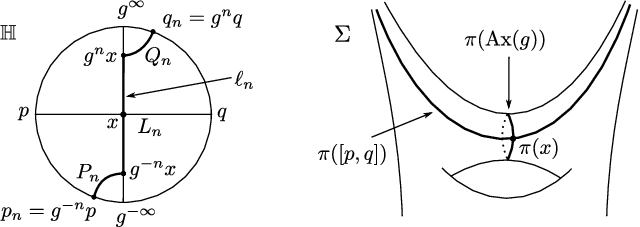}
\caption{In the right of the figure, we see that $\pi(\ell_n)=\pi([p,q])\cup \pi(\mathrm{Ax}(g))$ and the $L_n$-part of $\pi(\ell_n)$ goes round $\pi(\mathrm{Ax}(g))$ $2n$ times. Hence, it is intuitively natural that ``$\frac{1}{2n}\ell_n$'' converges to ``$g$''.}\label{fig:converging-path}
\end{figure}

Let $G_0$ be a subset of the complete system of representatives of $G/\langle g \rangle$ satisfying the condition that for any $h\in G_0$, the orbit $h\{ g^\infty , g^{-\infty}\}$ belongs to $A(B(K,r))$, i.e.,
\[ h\mathrm{Ax}(g)\cap B(K,r)\not=\emptyset.\]
Then we have
\[ \int f d\eta_g=\sum_{h\in G_0}f(h\{ g^{-\infty} , g^\infty\} ) \]
since $A(K)$ includes $\mathrm{supp}(f)$.

Now, we set
\[ H_0:=\{ h\in G\mid h[p,q]\cap B(K,r)\not=\emptyset \} \]
and
\[ G_n:=\{ h\in G\mid h\ell_n \cap B(K,r)\not=\emptyset \} \]
for $n\in \NN$.
For $h\in G$, if $h\{ p_n,q_n\} \in A(K)$, then $h[p_n,q_n]\cap K\not=\emptyset$ and hence $h\ell_n\cap B(K,r)\not=\emptyset $, which implies that $h\in G_n$. Therefore,
\[ \int f d\eta_{\{ p_n,q_n\} }=\sum_{h\in G_n }f(h\{ p_n ,q_n\} ).\]

Note that $h\ell_n\cap B(K,r)\not=\emptyset$ implies that
\[ hP_n\cap B(K,r)\not=\emptyset,\ hL_n\cap B(K,r)\not=\emptyset \text{ or }hQ_n\cap B(K,r)\not=\emptyset.\]
If $hP_n\cap B(K,r)\not=\emptyset$, then $hg^{-n}[p,q]\cap B(K,r)\not=\emptyset$ and so $h\in H_0g^n$.
Similarly, if $hQ_n\cap B(K,r)\not=\emptyset$, then $hg^n [p,q]\cap B(K,r)\not=\emptyset $ and so $h\in H_0g^{-n}$.
Now, we consider the case in which $hL_n\cap B(K,r)\not=\emptyset$.
Since $L_n=[g^{-n}x,g^nx]$ is included in $\mathrm{Ax}(g)$, we see that $h\mathrm{Ax}(g)\cap B(K,r)\not=\emptyset$, which implies that
\[ h \in \{ ug^m \mid u\in G_0, m\in \ZZ \}.\]
We remark that the number of $m\in \ZZ$ satisfying the condition that $ug^mL_n\cap B(K,r)\not=\emptyset$ is at most finite since $L_n$ is a finite geodesic segment.

We can assume that for every $u\in G_0$, we have $u[x,gx]\cap B(K,r)\not=\emptyset$.
Then for $u\in G_0$, we see that
\[ ug^mL_n\cap B(K,r)\not=\emptyset \text{ for } m\in [-n+1,n]\cap \ZZ .\]
Hence, for each $u\in G_0$, we can take a finite subset $Z_u^n$ of $\ZZ$ such that
\[ H_n:=\{ h\in G\mid hL_n\cap B(K,r)\not=\emptyset \} =\{ ug^m\mid u\in G_0, m\in Z_u^n\} ,\]
and $2n\leq \# Z_u^n \leq 2n+d$ for some constant $d>0$ depending on the diameter of $B(K,r)$.
Moreover, $Z_u^n$ includes $[-n+1,n]\cap \ZZ$.

From the above, we see that
\begin{align*}
&\left| \frac{1}{2n}\int f d \eta_{\{ p_n,q_n\} }-\int f d \eta_g \right|\\
=&\left| \frac{1}{2n}\sum_{h\in G_n}f (h\{ p_n,q_n\} )-\sum_{u\in G_0} f(u\{ g^{-\infty}, g^\infty \}) \right|\\
\leq &\left| \frac{1}{2n}\sum_{h\in H_n}f (h\{ p_n,q_n\} )-\sum_{u\in G_0} f(u\{ g^{-\infty}, g^\infty \} ) \right| \\
&\quad +\frac{1}{2n}\sum_{h\in H_0g^n\cup H_0g^{-n}}|f (h\{ p_n,q_n\} )|\\
\leq &\left| \sum_{u\in G_0}\sum_{m\in Z_u^n}\frac{1}{2n}f (ug^m\{ p_n,q_n\} )-\sum_{u\in G_0} f(u\{ g^{-\infty}, g^\infty \}) \right|
+\frac{2\# H_0}{2n}\max |f|\\
\leq &\sum_{u\in G_0}\left|\sum_{m\in Z_u^n}\frac{1}{2n}f (ug^m\{ p_n,q_n\} )-f(u\{ g^{-\infty}, g^\infty \}) \right| +\frac{\# H_0}{n}\max |f|\\
\end{align*}

Now, it is sufficient to see that for each $u\in G_0$, the sum
\[ \sum_{m\in Z_u^n}\frac{1}{2n}f (ug^m\{ p_n,q_n\} )=\sum_{m\in Z_u^n}\frac{1}{2n}f (u\{ g^{m-n}p,g^{m+n}q\} )\]
converges to $f(u\{ g^\infty ,g^{-\infty} \} )$ when $n$ tends to infinity.
The idea is almost the same as that of the proof of $\lim_{n\rightarrow \infty}\frac{1}{n}\sum_{k=1}^n a_k=\alpha$ for a sequence $\{ a_n\}$ of $\RR$ converging to $\alpha\in \RR$.
Fix $u\in G_0$ and $\varepsilon >0$. Since $f$ is continuous, there exists $N\in \NN$ such that if $j,k\geq N$, then
\[ | f(u\{ g^{-j}p, g^kq\} )-f(u\{ g^{-\infty },g^\infty \} )|<\varepsilon.\]
Note that we have the equation
\[ Y_u^n:=\{ m\in Z_u^n\mid m-n\leq -N \text{ and } m+n\geq N\} =[-n+N, n-N]\cap \ZZ \]
since $Z_u^n$ includes $[-n+1, n]\cap \ZZ$.
Therefore,
\begin{align*}
&\left| \sum_{m\in Z_u^n}\frac{1}{2n}f (ug^m\{ p_n,q_n\} )- f(u\{ g^\infty ,g^{-\infty} \} )\right|\\
\leq &\frac{1}{2n} \sum_{m\in Y_u^n}\big|f (u\{ g^{m-n}p,g^{m+n}q\} )-f(u\{ g^\infty ,g^{-\infty} \} )\big| \\
&\quad +\frac{1}{2n}\sum_{m\in Z_u^n\setminus Y_u^n}\Big( |f (u\{ g^{m-n}p,g^{m+n}q\} )|+|f(u\{ g^\infty ,g^{-\infty} \})| \Big)\\
< & \frac{2n-2N+1}{2n}\cdot \varepsilon+\frac{\#Z_u^n -\#Y_u^n }{2n}2\max |f|\\
\leq & \frac{2n-2N+1}{2n}\cdot \varepsilon+\frac{d+2N-1}{n}\max |f| \underset{n\rightarrow \infty }{\longrightarrow} \varepsilon.
\end{align*}
This completes the proof.
\end{proof}

\section{Discontinuity of intersection number on cusped hyperbolic surface}\label{sec:extension of intersection number}

For a compact hyperbolic surface $\gS$, the intersection number $i$ of closed geodesics was continuously extended to an $\RRR$-bilinear functional $i\: \GC (\gS )\times \GC(\gS) \rightarrow \RRR$
in \cite{Bon86}, i.e., for any closed geodesics $\gamma_1, \gamma_2$ on $\gS$, we have
\[ i(\eta_{\gamma_1},\eta_{\gamma_2})=i(\gamma_1,\gamma_2 ),\]
where $\eta_{\gamma_1}$ represents a counting geodesic current $\eta_c$ for $c\in \pi_1(\gS)$ satisfying the condition that a representative of $c$ is freely homotopic to $\gamma_1$. Note that if $c_1,c_2\in \pi_1(\gS)$ are conjugate, then $\eta_{c_1}=\eta_{c_2}$.

However, when $\gS$ is a cusped hyperbolic surface, we can see that the intersection number $i$ cannot extend continuously to an $\RRR$-bilinear functional $i\: \GC(\gS)\times \GC (\gS)\rightarrow \RRR$. The reason is as follows.
Let $\alpha, \beta \in \pi_1(\gS)$ be parabolic elements such that the fixed point $\alpha^\infty$ of $\alpha$ is different from the fixed point $\beta^\infty$ of $\beta$. From Proposition \ref{prop:converge seq} we have
\[ \eta_{\alpha^n\beta^n}\rightarrow 2\eta_{\{\alpha^\infty, \beta^\infty \} }\quad (n\rightarrow \infty ).\]
Then the intersection number of $\alpha^n\beta^n$ and the geodesic $\ell$ connecting two cusps corresponding to $\{\alpha^\infty ,\beta^\infty\}$ tends to infinity but the self-intersection number of $\ell$ is finite.

We remark that according to \cite[Theorem 2.4]{BIPP19}, for a cusped hyperbolic surface $\gS$ and any horocycle parameter $\lambda$
the intersection number
\[ i\: \GC (\gS )\times \GC_\lambda (\gS )\rightarrow \RRR \]
is continuous for
\[ \GC_\lambda:=\{ \mu\in \GC (\gS )\mid [ x,y]\subset \HH_\lambda \text{ for any }\{ x,y\} \in \mathrm{supp}(\mu) \}.\]
Roughly speaking, by restricting $\HH_\lambda$ we can prove the continuity of $i$ in the same manner as that in \cite{Bon86}.
In the rest of this section, we present another sketch of the proof of the continuity by using the argument in \cite[Section 5.3]{Sas19}, which was used for the proof of the continuity of the generalized intersection number on $\SC(\gS)$.

Let $\gS$ be a cusped hyperbolic surface. Fix a Dirichlet fundamental domain $\F$ for the action of $G=\pi_1(\gS)$ on $\HH$. By removing some edges of $\F$, we assume that
$G(\F)=\HH$ and $g\F\cap \F=\emptyset$ for any non-trivial $g\in G$. Set
\[ \mathcal{I}_\F=\{ (S_1,S_2)\in \partial_2\HH \times \partial_2\HH \mid CH(S_1)\cap CH(S_2)\text{ is a point in } \F\}.\]
Then we define
\[ i (\mu,\nu )=\mu\times \nu (\mathcal{I}_\F) \]
for any $\mu,\nu\in \GC (\gS)$. We see that $i(\eta_{\gamma_1},\eta_{\gamma_2})$ equals the intersection number of $\gamma_1$ and $\gamma_2$ for any closed geodesic or geodesic connecting two cusps $\gamma_1,\gamma_2$.
Remark that $\mu\times \nu (\mathcal{I}_\F)$ is independent of the choice of $\F$ since $\mu,\nu$ are $G$-invariant.

Fix any $(\mu,\nu)\in \GC (\gS )\times \GC_\lambda (\gS )$ and take a sequence $\{ (\mu_n,\nu_n)\}$ of $\GC (\gS )\times \GC_\lambda (\gS )$ converging to $(\mu, \nu)$.
From Portmanteau theorem \cite[Proposition 5.45]{Sas19}, if $\mu\times \nu (\partial\mathcal{I}_\F)=0$, then $\mu_n\times \nu_n(\mathcal{I}_\F)$ converges to $\mu\times \nu(\mathcal{I}_\F)$. However, in general, $\mu\times \nu (\partial\mathcal{I}_\F)$ is not necessarily zero.
Note that $(S_1,S_2)\in \partial\mathcal{I}_\F$ satisfies one of the following two conditions:
\begin{enumerate}
\item $S_1\not=S_2$ and $CH(S_1)\cap CH(S_2)$ is a point on $\partial \F$;
\item $S_1=S_2$ and $CH(S_1) \cap \ol{\F}\not=\emptyset$.
\end{enumerate}
From the assumption that $\nu\in \GC_\lambda$, for the condition (1), it is enough to consider the case in which $CH(S_1)\cap CH(S_2)$ is a point on $\partial \F\cap \HH_\lambda$; for the condition (2), it is enough to consider the case in which $CH(S_1) \cap (\ol{\F}\cap \HH_\lambda)\not=\emptyset$.

For the condition (1), by moving the center of $\F$ we can assume that the measure of such $(S_1,S_2)$ by $\mu \times \nu$ equals $0$ (see Lemma \cite[Lemma 5.51]{Sas19} for further details).

For the condition (2), we see that if the set of such $(S,S)$ has a non-zero measure for $\mu\times \nu$, then $\mu$ and $\nu$ has a common atom $\delta_S$ for $(S,S)\in \partial \mathcal{I}_\F$. Since $\nu \in \GC_\lambda$, $S$ is the limit set of $\langle g\rangle $ for a hyperbolic element $g\in G$.
Then we can prove that there exists a small open neighborhood $V$ of $(S,S)$ such that $\mu\times \nu (V\setminus \{ (S,S)\})$ is arbitrary small (see the proof of \cite[Theorem5.39]{Sas19} for further details).
As a result, we can prove that $\mu_n\times \nu_n (\mathcal{I}_\F)$ converges to $\mu\times \nu(\mathcal{I}_\F)$. This implies that the intersection number
\[ i\: \GC (\gS )\times \GC_\lambda (\gS )\rightarrow \RRR \]
is continuous.

\section{Proof of denseness property of rational geodesic currents}\label{sec:proof dense gc}

Let $\gS$ be a cusped hyperbolic surface (recall Assumption \ref{assump:gS has no boundary}).
The main purpose of this section is to prove that the space $\GC (\gS)$ has the denseness property of rational geodesic currents (Theorem \ref{thm:GC denseness 1}). From Theorem \ref{thm:dense of subgp and current}, it is sufficient to prove that the set $\GC^d(\gS)$ of discrete geodesic currents is a dense subset of $\GC (\gS)$.

Our strategy for the proof is based on the proof of the denseness property of rational geodesic currents on a hyperbolic group in \cite{Bon88} and the proof of the denseness property of rational subset currents on a surface group in \cite{Sas19}.
For a given geodesic current $\mu\in \GC( \gS )$, we construct a $G$-invariant family of quasi-geodesics on $\HH$, which induces a discrete geodesic current on $\gS$ approximating $\mu$ by considering the limit set of each quasi-geodesic.
To construct the quasi-geodesics, we introduce the notion of a ``round-path'', which is an analogy of a ``round-graph'' in \cite{Sas19}. Roughly speaking, by combining a $G$-invariant family of round-paths we can obtain a $G$-invariant family of quasi-geodesics on $\HH$.

Fix a fundamental domain $\F$ for the action of $G=\pi_1(\gS)$ on $\HH$ such that $\F$ is a convex polygon whose vertices on $\partial \HH$ are parabolic fixed points of $G$. We can obtain such a fundamental domain by cutting $\gS$ along some geodesics connecting two cusps.
We remove some edges of $\F$, which is a boundary component of $\F$, such that $G(\F)=\HH$ and $g\F\cap \F=\emptyset $ for any $g\in G\setminus \{ \id \}$.
In this setting, the set $\mathrm{Side}(\F)$ of side-pairing transformations of $\F$ is a basis of the free group $G$.
Then we consider the Cayley graph $\Cay(G)$ of $G$ with respect to the basis $\mathrm{Side}(\F)$.
Recall that the vertex set $V(\Cay (G))$ of $\Cay(G)$ is $G$; the edge set $E(\Cay (G))$ is $G\times \mathrm{Side}(\F)$, and an edge $(g,a)\in E(\Cay (G))$ connects $g$ to $ga$. We endow $\Cay(G)$ with the path metric $d_G$ such that every edge has length $1$.
Note that $g_1,g_2\in V(\Cay (G))$ are adjacent if and only if $g_1\F$ and $g_2\F$ are adjacent.

Let $B_G(g,r)$ be the closed ball centered at $g\in G=V(\Cay (G))$ with radius $r\geq 0$ in $\Cay (G)$.
Take some horocycle parameter $\lambda$ and set
\[ \F_\lambda :=\F \cap \HH_\lambda=\F\cap \pi^{-1}(\gS_\lambda ).\]
For $g\in G$ and $r\in \NN\cup \{ 0\}$, we consider the $r$-neighborhood
\[ B_G(g\F ,r)=\bigsqcup_{h\in V(B_G(g,r))}h\F \]
of $g\F$ with respect to $d_G$ and the $r$-neighborhood
\[ B_G(g\F_\lambda ,r)=B_G(g\F ,r)\cap \HH_\lambda =\bigsqcup_{h\in V(B_G(g,r))}h\F_\lambda\]
of $g\F_\lambda$ with respect to $d_G$.

For an edge $e$ of $\F$, which is included in $\ol{\F}$ but not necessarily included in $\F$, we call $e\cap \HH_\lambda$ an edge of $\F_\lambda$.
For a horocycle $U$ on $\HH$, which corresponds to a boundary component of $\gS_\lambda$, we call the intersection $U\cap \F$ a \ti{horocyclic edge} of $\F_\lambda$.
We say that $e$ is an edge of  $B_G(g\F_\lambda ,r)$ if $e$ is an edge of $h\F_\lambda$ for some $h\in B_G(g,r)$.

The notion of a round-path, which we define in the following, will play a fundamental role in the proof of the denseness property of rational geodesic currents.
Roughly speaking, a round-path of $B_G(g\F_\lambda ,r)$ is an information of how a geodesic line on $\HH$ passes through $B_G(g\F_\lambda ,r)$.

\begin{definition}[Round-path and Cylinder]
For a sequence of edges $e_1,\dots , e_k$ of $B_G(g\F_\lambda, r)$, we say that $[e_1,\dots ,e_k]$ is a \ti{round-path} of $B_G(g\F_\lambda ,r)$ if there exists a geodesic line $\ell$ in $\HH$ passing through $e_1,\dots , e_k$ in this order while passing through $\ol{B_G(g\F_\lambda,r)}$.
In this case we say that $\ell$ \ti{passes through} the round-path $[e_1, \dots ,e_k]$.
Note that we can also consider the case that a quasi-geodesic passes through $[e_1,\dots ,e_k]$ similarly.
We identify $[e_1,\dots ,e_k]$ with $[e_k,e_{k-1},\dots ,e_1]$ since we do not consider the direction of $\ell$.

We add a supplementary explanation to the above definition for the completeness.
When $\ell$ passes through the intersection point of an horocyclic edge $e$ and another non-horocyclic edge of $\F_\lambda$, we consider $\ell$ to be  passing through the horocyclic edge $e$. When $\ell$ passes through the intersection point $v$ of two horocyclic edges $e_1$ of $h_1\F_\lambda$ and $e_2$ of $h_2\F_\lambda$ and $v$ belongs to $h_1\F_\lambda$, we consider $\ell$ to be passing through $e_1$.
For a round-path $p=[e_1,\dots ,e_k]$ of $B_G(g\F_\lambda ,r)$, if $e_i$ is a horocyclic edge of $h_i\F_\lambda$, $e_{i+1}$ is a horocyclic edge of $h_{i+1}\F_\lambda$ and $h_1,h_2$ are neither the same nor adjacent in $\Cay (G)$, then we do not consider such a round-path $p$ and consider $[e_1,\dots ,e_i]$ or $[e_{i+1},\dots ,e_k]$ instead of $p$, i.e., if a geodesic line $\ell$ passes through $e_1,\dots ,e_k$ in this order while passing through $\ol{B_G(g\F_\lambda ,r)}$, then we consider $\ell$ to be passing through $e_1,\dots ,e_i$ or $e_{i+1},\dots ,e_k$.

For a round-path $p$ of $B_G(g\F_\lambda ,r)$, we define the \ti{cylinder} $\Cyl (p)$ with respect to $p$ to be the subset of $\partial_2\HH$ consisting of $S$ satisfying the condition that $CH(S)$ passes through $p$.
We denote by $\R_r(g)$ the set of all round-paths of $B_G(g\F_\lambda ,r)$ that contain an edge $e$ of $g\F_\lambda$ with $e\cap g\F_\lambda \not=\emptyset$. Note that $\R_r(g)$ is a finite set.

For a round-path $p=[e_1,\dots ,e_k]\in \R_r(g)$ and $h\in G$, we define $hp$ to be $[he_1,\dots ,he_k]\in \R_r(hg)$. Hence, we have the action of $G$ on the union $\bigsqcup_{g\in G} \R_r(g)$.
\end{definition}

\begin{example}
In Figure \ref{fig:round-path}, we present three examples of round-paths in the upper-half plane model of $\HH$.
In the left of Figure \ref{fig:round-path}, the geodesic line $\ell$ passes through $e_1,e_2,e_3$ and $e_4$ in this order, all of which are not horocyclic edges.
In the center of Figure \ref{fig:round-path}, the geodesic line $\ell$ passes through $e_1,e_2,e_3$ and $e_4$ in this order, and $e_2,e_3$ are horocyclic edges.
In the right of Figure \ref{fig:round-path}, the geodesic line $\ell$ passes through $e_1,e_2,e_2$ and $e_3$ in this order, and only $e_2$ is a horocyclic edge.
Note that a horocyclic edge is not a geodesic segment.

\begin{figure}[h]
\centering
\includegraphics{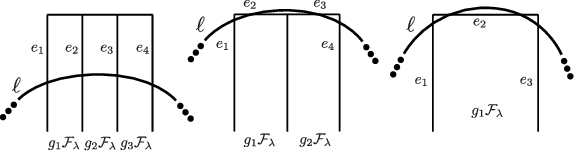}
\caption{Example of round-paths.}\label{fig:round-path}
\end{figure}

\end{example}

By the definition of round-paths, we see that for two different round-paths $p_1,p_2\in \R_r(g)$,
\[ \Cyl(p_1)\cap \Cyl(p_2)=\emptyset .\]
Moreover, we have
\[ \bigsqcup_{p\in \R_r(g)}\Cyl (p)=\{ S\in \partial_2\HH \mid CH(S)\cap g\F_\lambda \not=\emptyset \} =A(g\F_\lambda ).\]
since $CH(S)$ passes through an edge $e$ of $g\F_\lambda$ with $e\cap g\Fla \neq \emptyset$ for $S\in \partial_2\HH$ if and only if $CH(S)\cap g\F_\lambda\not=\emptyset$.

Now, we prepare some lemmas related to round-paths and cylinders.

\begin{lemma}\label{lem:distance between horo edges}
Take any $L>0$. If a horocycle parameter $\lambda$ is sufficiently large, then for any $r\in \NN, g\in G$ and $p=[e_1,\dots ,e_k]\in \R_r(g)$ such that $e_1,e_k$ are horocyclic edges, the distance
\[ d(e_1,e_k):=\inf \{ d(x,y)\mid x\in e_1, y\in e_k \}\]
is larger than $L$.
\end{lemma}
\begin{proof}
Fix some horocycle parameter $\lambda_0$. Assume that the horocycle parameter $\lambda$ is larger than $\lambda_0$ so that the distance between each boundary component of $\gS_{\lambda_0}$ and the boundary component of $\gS_{\lambda}$ corresponding to the same cusp is larger than $L/2$.
Then considering that $\gS_\lambda, \gS_{\lambda_0}$ are subsurfaces of $\gS$, we see that if a geodesic line $\ell$ on $\gS$ goes into $\gS_{\lambda}$ from a cusp neighborhood, then $\ell$ goes down the cusp neighborhood to $\gS_{\lambda_0}$.
Take any $x\in e_1$ and $y\in e_k$. Then $\pi ([x,y])$ goes into $\gS_{\lambda}$ at $x$ and goes out from $\gS_{\lambda}$ at $y$, which implies that $\pi ([x,y])$ passes through the cusp neighborhood between $\gS_{\lambda}$ and $\gS_{\lambda_0}$ twice. Therefore, the length of $[x,y]$ is larger than $L$.
\end{proof}

Let $d_{\ol{\HH}}$ be the distance function on $\ol{\HH}=\HH\cup \partial \HH$, which is the restriction of the Euclidean distance to the Poincar\'e disk model of $\HH$. We define the distance function $d_{\partial_2\HH}$ on $\partial_2\HH$ as
\[ d_{\partial_2\HH}(S_1,S_2)=\max \{ \max_{x\in S_1}d_{\ol{\HH}}(x,S_2), \max_{x\in S_2}d_{\ol{\HH}}(S_1,x)\}\]
for $S_1,S_2\in \partial_2\HH$, which is the restriction of the Hausdorff distance on $\H (\partial \HH)$ with respect to $d_{\ol{\HH}}$ to $\partial_2\HH$.

\begin{lemma}\label{lemma:diam of cyl}
Fix $g\in G$ and $\varepsilon >0$. If a horocycle parameter $\lambda$ and $r\in \NN$ is sufficiently large, then for any $p\in \R_r(g)$, the diameter of $\Cyl (p)$,
\[ \mathrm{diam}(\Cyl (p)):=\sup \{ d_{\partial_2\HH}(S_1,S_2)\mid S_1, S_2\in \Cyl (p)\},\]
is smaller than $\varepsilon$.
\end{lemma}
\begin{proof}
To obtain a contradiction, suppose that there exists $\varepsilon >0$ such that for any horocycle parameter $\lambda_0$ and $r_0\in \NN$, there exist $\lambda$ larger than
$\lambda_0$, $r\geq r_0$, $p\in \R_r(g)$ and $S_1,S_2\in \Cyl (p)$ such that $d_{\partial_2\HH}(S_1,S_2)\geq \varepsilon$.
We can assume that for some $x_1\in S_1$, we have $d_{\ol{\HH}}(x_1, S_2)>0$ without loss of generality. Then there exists $\varepsilon'>0$ depending only on $\varepsilon$ such that
$d_{\ol{\HH} }(x_1,CH(S_2))>\varepsilon'$. Hence, for the closed ball $B_{\ol{\HH}}(x,\varepsilon ')$ centered at $x$ with radius $\varepsilon'$, we have
\[ B_{\ol \HH}(x,\varepsilon ')\cap CH(S_2)=\emptyset.\]

Now, we assume that $r_0$ and $\lambda_0$ are sufficiently large such that for any $y\in \partial \HH$, some edges of $B_G(g\Fla, r)$ are included in $B_{\ol \HH}(y,\varepsilon')$.
Moreover, we can assume that an end edge $e$ of $p$ is included in $B_{\ol \HH} (x,\varepsilon ')$.
Since $S_2\in \Cyl (p)$, we have $CH(S_2)\cap e\not=\emptyset$. Hence, $CH(S_2)\cap B_{\ol \HH}(x, \varepsilon')\not=\emptyset$, which is a contradiction.
\end{proof}

From the proof of the above lemma and the stability of quasi-geodesics on a Gromov hyperbolic space, we can obtain the following lemma.
Recall that for $a\geq 1, b>0$, an $(a,b)$-quasi-geodesic on $\HH$ is an $(a,b)$-quasi-isometric embedding from an interval of $\RR$ to $\HH$. The quasi-geodesics that we are going to construct later are piecewise geodesics.

Recall that for a subset $A$ of $\HH$ the \ti{limit set} $A(\infty)$ of $A$ is the set of accumulation points of $A$ in $\partial \HH$.

\begin{lemma}\label{lemma:diam quasi-geod}
Fix $g\in G$, $a\geq 1, b\geq 0$ and $\delta_0>0$.
If a horocycle parameter $\lambda$ and $r\in \NN$ are sufficiently large, then for any $p=[e_1,\dots ,e_k]\in \R_r(g)$, if an $(a,b)$-quasi-geodesic $\ell$ passes through $e_1,\dots , e_k$ in this order while passing through $\ol{B_G(g\F_\lambda,r)}$, then
the limit set $\ell (\infty )$ of $\ell$ is contained in the $\delta_0$-neighborhood of $\Cyl(p)$.
\end{lemma}

\begin{proof}
To obtain a contradiction, assume that $\ell(\infty )$ is not contained in the $\delta_0$-neighborhood of $\Cyl(p)$.
Then we can take $\delta_0'>0$ depending on $\delta_0$, $S \in \Cyl(p)$ and $x\in S$ such that $B_{\ol{\HH}}(x,\delta_0 ')\cap CH(\ell (\infty ))=\emptyset$.
However, by the same argument as in the above lemma, if $\lambda$ and $r$ is sufficiently large, then $B_{\ol{\HH}}(x,\delta_0 ')\cap \ell \not=\emptyset$ and one of the limit points of $\ell$ is contained in $B_{\ol{\HH}}(x,\delta_0 ')$.
This implies that $B_{\ol{\HH}}(x,\delta_0 ')\cap CH(\ell (\infty ))\not=\emptyset$, which is a contradiction.
\end{proof}

Since the proof of the denseness property of rational geodesic currents is long and includes many constants and parameters, we will write
\textbf{Setting} when we fix something related to the proof.

\begin{setting8}
Fix $\mu\in \GC (\gS) $ and assume that $\mu$ is not the zero measure.
Fix $\varepsilon >0$. Take any continuous functions $f_1,\dots ,f_l\: \partial_2\HH\rightarrow \RRR$ with compact supports. Take the neighborhood of $\mu$ as follows:
\[ U(\varepsilon ; f_1,\dots ,f_l)=\left\{ \nu \in \GC (\gS )\mathrel{}\middle|\mathrel{} \left| \int f_i d\mu -\int f_i d\nu \right|<\varepsilon \ (i=1,\dots ,l)\right\}.\]
Take a compact subset $K$ of $\HH$ such that
\[ A(K)=\{ S\in \partial_2\HH \mid CH(S)\cap K\not=\emptyset \} \]
includes the support $\mathrm{supp}(f_i)$ of $f_i$ for $i=1,\dots ,l$.
From now on we assume that the horocycle parameter $\lambda$ is large enough so that $K$ is included in $\HH_\lambda$. In addition, we take $r_0\in \NN$ such that $K$ is included in
\[ B_G(\F_\lambda ,r_0)=\bigsqcup_{h\in V(B_G(\id ,r_0))}h\F_\lambda. \]

\end{setting8}
Note that the family of $U(\varepsilon ; f_1,\dots ,f_l)$ forms a fundamental system of neighborhoods of $\mu$.
We are going to construct a discrete geodesic current, i.e., a finite sum of rational geodesic currents
\[ \nu= c_1\eta_1 +\cdots +c_t\eta_t \ (c_1,\dots ,c_t>0),\]
belonging to the neighborhood $U(\varepsilon ; f_1,\dots ,f_l)$.

\begin{lemma}\label{lem:disjoint union of cylinder}
There exists a subset $\mathcal{O}$ of
\[ \bigsqcup_{h\in V(B_G(\id ,r_0 ))}\R_r (h)\]
such that
\begin{equation}
\bigcup_{h\in V(B_G(\id, r_0))} A(h\Fla )=\bigsqcup_{p\in \mathcal{O}}\Cyl (p). \tag{$\ast $}
\end{equation}
\end{lemma}
\begin{proof}
Recall the discussion before Lemma \ref{lem:distance between horo edges}. For each $h\in V(B_G(\id ,r_0 ))$, we have
\[ A(h\Fla )=\bigsqcup_{p\in \R_r (h)}\Cyl (p).\]
First, we set
\[ \mathcal{O}=\bigsqcup_{h\in V(B_G(\id ,r_0 ))}\R_r (h),\]
and we remove some round-paths from $\mathcal{O}$ such that $\mathcal{O}$ satisfies the above condition $(\ast)$.

Taking a labeling of the elements of $V(B_G(\id ,r_0 ))$, we have
\[ V(B_G(\id ,r_0 ))=\{ g_1, \dots ,g_s\} .\]
For $p_1\in \R_r (g_{i_1}) ,p_2\in \R_r(g_{i_2} )$, if $\Cyl (p_1)\cap \Cyl (p_2)\not=\emptyset$ and $i_1<i_2$, then we remove $p_2$ from $\mathcal{O}$.
We continue this operation for each pair of $p_1,p_2\in \mathcal{O}$ one by one. Finally, we can obtain $\mathcal{O}$ such that for any $p_1,p_2\in \mathcal{O}$ so that $p_1\not= p_2$, we have $\Cyl (p_1)\cap \Cyl (p_2)=\emptyset$.

Now, it is sufficient to prove that
\[ \bigcup_{h\in V(B_G(\id, r_0))} A(h\Fla )\subset \bigsqcup_{p\in \mathcal{O}}\Cyl (p).\]
Take any $S\in \bigcup A(h\Fla )$. Let $i$ be the smallest number in $\{ 1,\dots, s\}$ such that $CH(S)$ passes through $g_i \Fla$. Then we can take $p\in \R_r(g_i)$ such that $S\in \Cyl (p)$. Since $CH(S)$ does not pass through $g_1\Fla, \dots ,g_{i-1}\Fla$, the round-path $p$ does not have an edge $e$ intersecting
\[ g_1\Fla\sqcup \cdots \sqcup g_{i-1}\Fla.\]
Therefore, for any $p'\in \R_r(g_1)\sqcup \cdots \sqcup \R_r(g_{i-1} )$, we have $\Cyl (p)\cap \Cyl (p')=\emptyset$, which implies that $p$ is not removed from the original $\mathcal{O}$ in the above operation. Hence, $p\in \mathcal O$, and
\[ S\in \Cyl (p)\subset \bigsqcup_{p'\in \mathcal{O}}\Cyl (p'),\]
as required.
\end{proof}

\begin{notation}
Let $m$ be a Borel measure on a topological space $\Omega$. Set $|m| := m(\Omega)$. For a non-empty
Borel subset $A$ of $\Omega$, we denote by $m|_A$ the restriction of $m$ to $A$, i.e., for any Borel subset $E$ of $\Omega$, 
\[ m|_A(E):=m(A\cap E).\]
The support of $m$, denoted by $\mathrm{supp}(m)$, is the smallest closed subset $A$ of $\Omega$ such that $m(A^c) =0$.
\end{notation}

The following lemma will play a fundamental role in proving that a certain geodesic current $\nu$ belongs to the neighborhood $U(\varepsilon ; f_1,\dots ,f_l)$ of $\mu$.

\begin{lemma}\label{lem:approximation data}
There exist a horocycle parameter $\lambda$, a radius $r\in \NN$ of round-path, $\hat \varepsilon>0$ and $\delta_0>0$ such that
if a geodesic current $\nu\in \GC (\gS)$ satisfies the following conditions, then $\nu$ belongs to $U(\varepsilon ; f_1,\dots ,f_l)$:
\begin{enumerate}
\item Take $\mathcal O$ satisfying the condition $(\ast)$ in Lemma \ref{lem:disjoint union of cylinder}. There exists a Borel measure $\nu_p$ for each $p\in \mathcal O$ such that
\[ \nu |_{A(K)}=\sum_{p\in \mathcal O}\nu_p |_{A(K)};\]
\item $\mathrm{supp}(\nu_p) $ is included in  the $\delta_0$-neighborhood $B(\Cyl (p),\delta_0 )$ of $\Cyl (p)$ for every $p\in \mathcal O$;
\item $\big| |\nu_p|-\mu (\Cyl (p) )\big| <\hat \varepsilon$ for every $p\in \mathcal O$.
\end{enumerate}
\end{lemma}
\begin{proof}
Let $f\in \{ f_1,\dots ,f_l\}$. Take $\mathcal O$ satisfying the condition $(\ast)$ in Lemma \ref{lem:disjoint union of cylinder}.
Recall that
\[ \mathrm{supp}(f)\subset A(K)\subset \bigcup_{h\in V(B_G(\id ,r_0))}A(h\Fla )=\bigsqcup_{p\in \mathcal O}\Cyl (p).\]
Hence, we have
\begin{align*}
\left| \int fd\nu-\int fd\mu \right|
=&\left| \int f d\sum_{p\in \mathcal O}\nu_p -\sum_{p\in \mathcal O}\int_{\Cyl (p)}f d\mu\right| \\
\leq &\sum_{{\substack{p\in \mathcal{O}\\[1pt] \Cyl(p)\cap A(K)\not=\emptyset }}}\left| \int fd\nu_p -\int_{\Cyl (p)}fd\mu \right|.
\end{align*}

Since $f$ is uniformly continuous, for $\varepsilon_2>0$, there exists $\varepsilon_1>0$ such that
\[ \sup_{{\substack{x,y\in \partial_2\HH\\[1pt] d(x,y)<\varepsilon_1}}}| f(x)-f(y)| <\varepsilon_2.\]
Take $0<\delta_0 <\varepsilon_1$. From Lemma \ref{lemma:diam of cyl} there exist $\lambda$ and $r$ such that the diameter $\mathrm{diam}(\Cyl(p))$ is smaller than $\varepsilon_1-\delta_0$ for any $p\in \mathcal{O}$, which implies that the diameter of $B(\Cyl (p), \delta_0)$ is smaller than $\varepsilon_1$.

For each $p\in \mathcal O$, take some $x_p\in \Cyl (p)$. Then we have
\begin{align*}
&\left| \int fd\nu_p -\int_{\Cyl (p)}fd\mu \right|\\
\leq &\left| \int fd\nu_p -f(x_p)|\nu_p |\right| +\left| f(x_p)|\nu_p |-\int_{\Cyl (p)}fd\mu \right|\\
\leq &\varepsilon_2 |\nu_p |+\big| f(x_p)|\nu_p |-f(x_p)\mu (\Cyl (p))\big| +\left| f(x_p)\mu(\Cyl (p))-\int_{\Cyl(p)} fd\mu\right|\\
\leq & \varepsilon_2 |\nu_p |+|f(x_p)| \cdot \big| |\nu_p|-\mu (\Cyl (p) )\big| +\varepsilon_2\mu (\Cyl (p)).
\end{align*}
Set $\mathcal O':= \{ p\in \mathcal O \mid \Cyl (p)\cap A(K)\not=\emptyset \}$.
We obtain
\begin{align*}
&\sum_{p\in \mathcal O'}\left| \int fd\nu_p -\int_{\Cyl (p)}fd\mu \right| \\
\leq &\varepsilon_2\sum_{p\in \mathcal O'} (|\nu_p |+\mu (\Cyl (p)) )
	+\mathrm{max}f\sum_{p\in \mathcal O'} \big| |\nu_p|-\mu (\Cyl (p) )\big| \\
<& \varepsilon_2 \sum_{p\in \mathcal O'} (2\mu (\Cyl (p) )+\hat \varepsilon )+\mathrm{max}|f|\cdot \# \mathcal O' \cdot \hat \varepsilon\\
<& 2\varepsilon_2 \mu (B(A(K),\varepsilon_1 ))+(\varepsilon_2 +\mathrm{max}|f| )\hat \varepsilon \cdot \# \mathcal O'
\end{align*}
Note that $\Cyl (p)$ is included in the $\varepsilon_1$-neighborhood $B(A(K),\varepsilon_1)$ of $A(K)$ for $p\in \mathcal{O}'$.
Since $\mu$ is a regular measure, the value $\mu (B(A(K),\varepsilon_1 ))$ is close to $\mu (A(K))$ when $\varepsilon_1$ is small. Hence, we can consider $\mu (B(A(K),\varepsilon_1 ))$ as given. In addition, if we take a sufficiently small $\varepsilon_2$, which influences $\lambda$ and $r$, then $\varepsilon_2 \mu (B(A(K),\varepsilon_1 )) $ is smaller than $\varepsilon/2$.
Remark that the cardinality $\# \mathcal O'$ can become larger when $\lambda$ and $r$ become larger.
Therefore, we take sufficiently small $\hat \varepsilon$ after fixing $\lambda$ and $r$.
As a result, we can see that
\[ \left| \int fd\nu-\int fd\mu \right|< \frac{\varepsilon}{2} +\frac{\varepsilon}{2}< \varepsilon,\]
as required.
\end{proof}

\begin{setting8}\label{setting2}
Fix a horocycle parameter $\lambda$, a radius $r\in \NN$ of round-path, $\hat \varepsilon>0$ and $\delta_0>0$ as in the above Lemma.
From Lemma \ref{lemma:diam quasi-geod} and the proof of the above lemma, for some $a\geq 1, b\geq 0$, we can also assume that
for any $p=[e_1,\dots ,e_k]\in \mathcal{O}$ with $p\in \R_r(g)$, if an $(a,b)$-quasi-geodesic $\ell$ passes through $e_1,\dots , e_k$ in this order while passing through $\ol{B_G(g\F_\lambda,r)}$, then the limit set $\ell (\infty )$ of $\ell$ is included in the $\delta_0$-neighborhood $B(\Cyl (p),\delta_0 )$ of $\Cyl(p)$.
Note that the number of $g\in G$ satisfying the condition that $p\in \R_r(g)$ for some $p\in \mathcal{O}$ is finite.
\end{setting8}

\begin{definition}[Connectability]
Let $u,v\in G$ that are adjacent in $\Cay (G)$.
Let $p=[e_1,\dots ,e_k]\in \R_r(u)$ and assume that $p$ passes through an edge $e$ of $v\Fla$ with $e\cap v\Fla\not=\emptyset$.
Then the restriction of $p$ to 
\[ B_G(u\Fla, v\Fla, r): = B_G(u\Fla ,r)\cap B_G(v\Fla ,r),\]
denoted by $p|_{u,v}$, is defined as a sub-round-path $[e_i, e_{i+1},\dots ,e_j]$ of $p$ if the edges of $p$ included in $\ol{B_G(u\Fla ,v\Fla, r)}$ are $e_i,e_{i+1},\dots ,e_j$.
We call $[e_i, e_{i+1}\dots ,e_j]$ a \ti{round-path} of $B_G(u\Fla ,v\Fla, r)$.
We remark that a round-path of $B_G(u\Fla ,v\Fla, r)$ always includes an edge $e_u$ of $u\Fla$ with $e_u\cap u\Fla \neq \emptyset$ and an edge $e_v$ of $v\Fla$ with $e_v\cap v\Fla \neq \emptyset$.
In addition, there exists $p'\in \R_r(v)$ such that the restriction $p'|_{u,v}$ of $p'$ to $B_G(u\Fla, v\Fla, r)$ equals $p|_{u,v}$.

For $p_1\in \R_r(u)$ and $p_2\in \R_r(v)$, we say that $p_1$ and $p_2$ are \ti{connectable} if $p_1$ includes an edge $e_v$ of $v\Fla$ with $e_v\cap v\Fla \neq \emptyset$, $p_2$ includes an edge $e_u$ of $u\Fla$ with $e_u\cap u\Fla \neq \emptyset$, and
$p_1|_{u,v}=p_2|_{u,v}$.
\end{definition}

For $\mu \in \GC(\gS)$, we define the map
\[ \ol{\mu}\: \bigsqcup_{g\in G}\R_r(g)\rightarrow \RRR\]
as
\[ \ol{\mu} (p):= \mu (\Cyl (p) ) \text{ for } p\in \bigsqcup_{g\in G}\R_r(g) .\]
Note that since $\mu$ is $G$-invariant, the map $\ol{\mu}$ is determined by a finite number of the values $\{ \ol{\mu} (p) \}_{p\in \R_r (\id )}$.

For adjacent $u,v\in G$ and any round-path $J$ of $B_G(u\Fla ,v\Fla, r)$, we have
\[ \bigsqcup_{\substack{p\in \R_r(u)\\[1pt] p|_{u,v}=J}}\Cyl (p)=\bigsqcup_{\substack{p'\in \R_r(v)\\[1pt] p'|_{u,v}=J}}\Cyl (p')\]
since each side of the equation can be considered as the cylinder with respect to $J$.
Hence, we can obtain
\[
\sum_{\substack{p\in \R_r(u)\\[1pt] p|_{u,v}=J}}\mu( \Cyl (p))=\sum_{\substack{p'\in \R_r(v)\\[1pt] p'|_{u,v}=J}}\mu(\Cyl ( p' ))
\]
and hence
\begin{equation}
\sum_{\substack{p\in \R_r(u)\\[1pt] p|_{u,v}=J}}\ol{\mu} (p)=\sum_{\substack{p'\in \R_r(v)\\[1pt] p'|_{u,v}=J}}\ol{\mu}(p').  \tag{$\ast_J$}
\end{equation}

By considering the action of $G$ on the Cayley graph $\Cay (G)$, the system of the equations $(\ast_J )$ for all adjacent $u,v\in G$ and all round-path $J$ of $B_G(u\Fla ,v\Fla, r)$ can be
considered as a finite homogeneous system of linear equations with respect to the variables $\ol{\mu} (p)$ for $p\in \R_r (\id )$.
Since the coefficients of these equations are integer, by \cite[Lemma 8.11]{Sas19} there exists a rational solution approximating $\ol{\mu}$, which induces a map
\[ \theta \: \bigsqcup_{g\in G}\R_r(g)\rightarrow \mathbb{Z}_{\geq 0}\]
satisfying the following conditions:
\begin{enumerate}
\item $\theta $ is $G$-invariant, i.e., for any $p\in \bigsqcup_{g\in G}\R_r(g)$ and $h\in G$, we have $\theta(p)=\theta(hp)$;
\item there exists $M\in \NN$ such that for any $p\in \bigsqcup_{g\in G}\R_r(g)$, we have
\[ \left| \frac{1}{M}\theta (p)-\mu (\Cyl (p))\right| < \hat{\varepsilon};\]
\item for any adjacent $u,v\in G$ and any round-graph $J$ of $B_G(u\Fla ,v\Fla, r)$,
\begin{equation}
\sum_{\substack{p\in \R_r(u)\\[1pt] p|_{u,v}=J}}\theta (p)=\sum_{\substack{p'\in \R_r(v)\\[1pt] p'|_{u,v}=J}}\theta (p').  \tag{$\ast_J'$}
\end{equation}
\end{enumerate}
The point is that $\frac{1}{M}\theta$ approximates $\ol{\mu}$ and satisfies the same property.

Now, considering $\theta(p)$ copies of round-paths $p$ for each $p\in \bigsqcup_{g\in G}\R_r(g)$, we will construct a family of $G$-invariant quasi-geodesics by combining the round-paths modulo the Equation $(\ast_J')$, which will induce a discrete geodesic current $\eta_\gG$.
Then
\[  \nu=\frac{1}{M}\eta_\gG \]
will satisfy the condition in Lemma \ref{lem:approximation data}, which implies that $\nu$ belongs the neighborhood $U(\varepsilon ; f_1,\dots ,f_l)$.

First, in the same manner as in \cite[Theorem 8.12]{Sas19}, we construct a graph $\gG$ that $G$ acts on.
We define the vertex set $V(\gG )$ of $\gG$ to be the set
\[ \{ v(g,p,i)\} _{g\in G, p\in \R_r(g), i=1,\dots ,\theta(p)} .\]
We regard $v(g,p,i)$ as a copy of $v(g,p,1)$ for $i=2,\dots, \theta(p)$ and we write it $v(g,p)$ for short when no confusion arises. When $\theta(p)=0$, there exists no vertex $v(g,p,i)$.
Define an action of $G$ on $\gG$ as
\[ hv(g,T,i)=v(hg,hT,i)\]
for $h\in G$ and $v(g,T,i)\in V(\gG)$.
Define a map $\iota$ from $V(\gG )$ to $V(\Cay (G))$ to be the natural projection, i.e., for $v(g,t,i)\in V(\gG)$,
\[ \iota (v(g,t,i))=g.\]

We define the edge set $E(\gG)$ by connecting two vertices in $V(\gG)$ $G$-equivariantly in the following way.
For each $u\in \mathrm{Side}(\F)$ and each round-path $J$ of $B_G(\Fla, u\Fla ,r)$, we connect a vertex $v(\id ,p,i)$ to a vertex $v(u,p',i')$ such that
\[ p|_{u,v}=J=p'|_{u,v}.\]
Since for each round-path $J$ of $B_G(\Fla, u\Fla ,r)$, the Equation $(\ast_J')$ holds, the number of vertices $v(\id , p, i )\in \iota^{-1}(\id )$ with $p|_{u,v}=J$ is equal to the number of vertices $v(u,p',i')\in \iota^{-1}(u)$ with $p'|_{u,v}=J$.
Hence, there exists a one-to-one correspondence between
\[ \{ v(\id ,p,i)\in \iota^{-1}(\id )\mid p|_{u,v}=J\} \text{ and } \{ v(u ,p' ,i')\in \iota^{-1}(u )\mid p'|_{u,v}=J\}.\]
Then we spread the above edges by the action of $G$.
Explicitly, for two vertices $v(\id, p,i)$ and $v(u,p',i')$ connected by an edge, we connect $hv(\id,p,i)$ to $hv(u,p',i')$ by an edge for every $h\in G$.

From the above, we obtain a graph $\gG$ that $G$ acts on; moreover, the $G$-equivariant map $\iota\: V(\gG)\rightarrow V(\Cay (G))$ naturally extends to the $G$-equivariant map
$ \iota\: \gG \rightarrow \Cay (G)$
satisfying the condition that the restriction of $\iota$ to each connected component of $\gG$ is injective.
From the above construction, for two adjacent $u,v\in V(\Cay (G))$, if two vertex $v(u,p)$ and $v(v,p')$ are connected by an edge, then $p$ and $p'$ are connectable.

\begin{lemma}\label{lem:Y is a segment or line}
For every vertex $v\in V(\gG)$, the degree of $v$ is smaller than or equal to $2$, which implies that a connected component $Y$ of $\gG$ is a point or homeomorphic to an interval of $\RR$.
Moreover, if a connected component $Y$ of $\gG$ is not a finite subgraph, then $Y$ is not a half-line but a line, i.e., homeomorphic to $\RR$.

In addition, if $v(g,p)\in V(\gG)$ is an end vertex of a finite connected component $Y$ of $\gG$, $p=[e_1,\dots ,e_k]$ and $e_1$ is an edge of $g\Fla$, then $e_1$ is a horocyclic edge of $g\Fla$.
\end{lemma}
\begin{proof}
Let $v(g, p)\in V(\gG)$. For the round-path $p\in \R_r(g)$, there exists a geodesic line $\ell$ passes through $p$.
If the vertex $v(g,p)$ is connected to $v(gu,p)$ for $u\in \mathrm{Side}(\F)\sqcup \mathrm{Side}(\F)$, then $\ell$ passes through $gu\Fla$.
By the definition of the fundamental domain $\F$, the number of such $u$ are at most two.
Hence, the degree of $v(g,p)$ is smaller than or equal to $2$.

Next, we consider a connected component $Y$ of $\gG$ that is not a finite subgraph.
The point is that $G$ acts on the graph $\gG$ and the quotient graph $G\backslash \gG$ is a finite graph since we have
\[ \# V(G\backslash \gG )=\# \iota^{-1}(\id )=\sum_{p\in \R_r(\id )}\theta (p)<\infty.\]
Hence, the quotient graph of $Y$ by the stabilizer of $Y$ with respect to the action of $G$ is also a finite graph, which implies that $Y$ can not be a half-line.
\end{proof}

Consider a connected component $Y$ of $\gG$ and assume that $Y$ has infinite vertices.
Since $Y$ is a line, we can assign a number to the vertex set $V(Y)$ of $Y$ such that
\[ V(Y)=\{ v(g_i,p_i)\} _{i\in \ZZ}\]
and $v(g_{i\pm 1},p_{i\pm 1})$ is connected to $v(g_i,p_i)$ for any $i\in \ZZ$.
Moreover, we can obtain a bi-infinite sequence $[ e_i]_{i\in \ZZ}$ of edges by combining the round-paths $\{ p_i\}_{i\in \ZZ}$ since adjacent round-paths of $\{ p_i\}_{i\in \ZZ}$ are connectable.

Even when $Y$ has at most finitely many vertices, we can obtain a finite sequence $\{ v(g_i,p_i)\}$ of vertices and a finite sequence $[e_i]$ of edges in the same manner.

\begin{lemma}\label{lem:combine round-path}
Let $Y$ be a connected component of $\gG$. There exists an infinite piecewise geodesic $\ell(Y)$ passing through the sequence $[e_i]$ of edges in this order
such that every bending angle of $\ell(Y)$ is larger than $\pi/2$ and every geodesic piece of $\ell(Y)$ is long enough that $\ell(Y)$ is an $(a,b)$-quasi-geodesic line $($see Setting \ref{setting2} for the constants $a,b)$.
\end{lemma}
\begin{proof}
1) First, we consider the case in which $Y$ is a line.
We assume that every $e_i$ is not a horocyclic edge for convenience.
For a finite subsequence $[e_{i-r},\dots ,e_{i+r}]$ of $[e_i]$, there exists a round-path $p_j$ such that $[e_{i-r},\dots ,e_{i+r}]$ is a subsequence of $p_j$. Hence, there exists a geodesic segment $\ell_i$ starting from $e_{i-r}$ toward $e_{i+r}$ and passing through $e_{i-r},\dots ,e_{i+r}$ in this order.
We combine the sequence $\{\ell_{ir}\}_{i\in \ZZ}$ of geodesic segments in the following way (see Figure \ref{fig:combine-round-path}) and construct an infinite piecewise geodesic $\ell(Y)$ passing through the sequence $[e_i]$ of edges in this order. Assume that $r$ is a multiple of $4$ for convenience.

\begin{enumerate}
\item[(1)] If $\ell_{ir}$ and $\ell_{ir+r}$ intersect at $t$ while passing through $e_{ir+\frac{r}{4}}$ and $e_{ir+\frac{3r}{4}}$, then
we combine $\ell_{ir}$ and $\ell_{ir+r}$ at $t$.
\item[(2)] If $\ell_{ir}$ and $\ell_{ir+r}$ do not intersect while passing through $e_{ir+\frac{r}{4}}$ and $e_{ir+\frac{3r}{4}}$, then we take the intersection point $s$ of $\ell_{ir}$ and $e_{ir+\frac{r}{4}}$ and take the intersection point $t$ of $\ell_{ir+r}$ and $e_{ir+\frac{3r}{4}}$, and we combine $\ell_{ir}$ with the geodesic segment $[s,t]$ at $s$ and combine $[s,t]$ with $\ell_{ir+r}$ at $t$.
\end{enumerate}

\begin{figure}[h]
\centering
\includegraphics{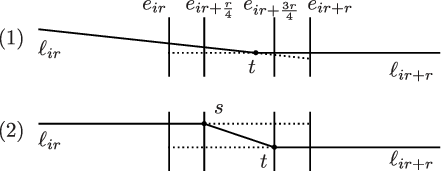}
\caption{The way of combining geodesic segments $\{ \ell_{ir}\}_{i\in \ZZ}$.}\label{fig:combine-round-path}
\end{figure}

From the above, we see that the length of every geodesic piece of $\ell(Y)$ is larger than
\[ \frac{r}{2}\inf \{ d(x,y)\mid x,y \text{ belong to non-adjacent edges of }\Fla \} \]
Each bending angle of $\ell(Y)$ is larger than $\pi/2$ if $r$ is sufficiently large.
Hence, if $r$ is sufficiently large, then $\ell(Y)$ is an $(a,b)$-quasi-geodesic line (see Supplementation \ref{supp:piecewise geodesic is quasi-geodesic}).

2) Now, we consider the case in which $Y$ is a finite segment. Let $[e_0,\dots ,e_m]$ be the finite sequence of edges that we obtained by combining the round-paths.
Then the end edges $e_0$ and $e_m$ must be horocyclic edges by the construction of $\gG$.
If $m$ is smaller than or equal to $2r$, then there exists a geodesic segment $[s,t]$ starting from $e_0$ toward $e_m$ and passing through $e_0,\dots ,e_m$ in this order. Then we combine $[s,t]$ with the geodesic ray $[s ,\xi]$ starting from $s$ to the parabolic fixed point $\xi$ that the horocycle including $e_0$ is centered at. Similarly, we combine $[s,t]$ with the geodesic ray $[t,\zeta ]$ starting from $t$ to the parabolic fixed point $\zeta$ that the horocycle including $e_m$ is centered at (see Figure \ref{fig:combine-round-path2}). This piecewise geodesic $\ell(Y)$ satisfies the condition in the lemma since the length of $[s,t]$ is sufficiently large by Lemma \ref{lem:distance between horo edges}.

\begin{figure}[h]
\centering
\includegraphics{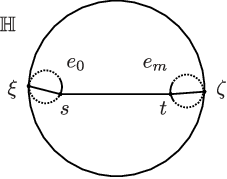}
\caption{The way of extending $[s,t]$ to $\ell(Y)$.}\label{fig:combine-round-path2}
\end{figure}

If $m$ is larger than $2r$, then we take geodesic segments
\[ \ell_r,\ell_{2r},\dots ,\ell_{cr}\quad (0\leq m-cr<r)\]
and combine them in the same manner as the above. Note that $\ell_{cr}$ is a geodesic segment starting from $e_{cr-r}$ to $e_m$.
Then we combine the resulting piecewise geodesic with the geodesic rays from its ends to the corresponding parabolic fixed points in the same manner as the above.
The resulting piecewise geodesic $\ell(Y)$ satisfies the condition in the lemma if $r$ is sufficiently large.
\end{proof}

\begin{supplementation}\label{supp:piecewise geodesic is quasi-geodesic}
It is well-known that a piecewise geodesic whose every bending angle is bounded from below and each segment is sufficiently long is a quasi-geodesic but we could not find any literature on this claim. Hence, we give the proof here for the convenience of the reader.

We use the fact that a local quasi-geodesic is a quasi-geodesic (see \cite[p. 25]{CDP90}).
Since each segment of the piecewise geodesic $\ell$ can be sufficiently long, it is enough to see that the neighborhood of each corner is a quasi-geodesic.
Therefore, we consider the case of Figure \ref{fig:piecewise-geodesic-evaluation}.
By the trigonometry of right triangles in the hyperbolic plane, we have
\[ \cosh (d(v,x))=\frac{\cos \alpha_1}{\sin \beta_1}=\frac{\cos \alpha_2}{\sin \beta_2}.\]
Therefore,
\[ \cosh (d(v,x))\leq \frac{1}{\max \{ \sin \beta_1, \sin \beta_2\} } \leq \frac{1}{\sin \frac{\beta_1 +\beta_2}{2} }.\]
This implies that $d(v,x)$ is bounded above by a constant $C$ depending on the bending angle $\beta_1+\beta_2$.
Then we see that
\[ d(u,v)+d(v,w)\leq d(u,x)+d(x,v) + d(v,x)+d(x,w) \leq d(u,w)+2C.\]
Hence, the neighborhood of the corner $v$ of $\ell$ is a $(1,2C)$-quasi-geodesic. This completes the proof.

\begin{figure}[h]
\centering
\includegraphics{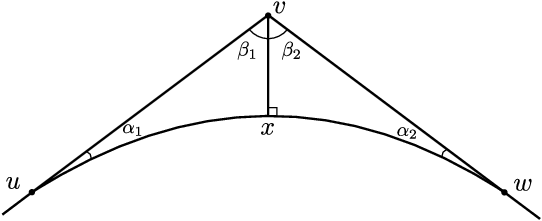}
\caption{The vertex $v$ of the triangle is a corner of $\ell$ and the vertices $u,w$ are on the segments of $\ell$.}\label{fig:piecewise-geodesic-evaluation}
\end{figure}

\end{supplementation}

Define $\Comp (\gG)$ to be the set of all connected components of $\gG$. Note that the action of $G$ on $\gG$ induces the action of $G$ on $\Comp (\gG)$.
However, since for each $Y\in \Comp (\gG)$, the construction of $\ell(Y)$ is not $G$-equivariant in the above proof, we add some supplementary explanation.

For $Y\in \Comp(\gG)$ homeomorphic to $\RR$ and an infinite quasi-geodesic $\ell (Y)$ satisfying the condition in the above lemma, 
the limit set $\ell(Y)(\infty )\in \partial_2\HH$ is uniquely determined by the bi-infinite sequence of the round-paths ${p_i}$. Hence, we see that $g\ell(Y)(\infty)$ equals $\ell(gY)(\infty)$ for any $g\in G$.
When $Y$ is finite, the limit set $\ell(Y)(\infty)$ is determined by the end edges of $[ e_i]$, which implies that $g\ell(Y)(\infty)$ equals $\ell(gY)(\infty)$ for any $g\in G$

Therefore, we can obtain a family $\{ \ell (Y)\}_{Y\in \Comp (\gG)}$ of infinite quasi-geodesics such that $\ell (Y)$ satisfies the condition in the above lemma for every $Y\in \Comp (\gG)$ and $G$ acts on $\{ \ell (Y)(\infty)\}_{Y\in \Comp (\gG)}$, i.e., for any $Y\in \Comp (\gG )$ and any $g\in G$
\[ \ell (gY)(\infty )=g\ell (Y)(\infty ).\]
As a result, we can obtain the $G$-invariant measure
\[ \eta_\gG =\sum_{Y\in \Comp (\gG )} \delta_{\ell (Y)(\infty )}\]
on $\partial_2 \HH$.
We will see that $\nu:=\frac{1}{M} \eta_\gG$ satisfies the condition in Lemma \ref{lem:approximation data}.

(Condition (1) and (2) in Lemma \ref{lem:approximation data}) For each $p\in \mathcal{O}$ with $p\in \R_r(g)$, we set
\[ \eta_{\gG,p}=\sum_{v(g,p)\in V(Y)} \delta_{\ell(Y)(\infty)},\]
where the sum is taken over all $Y\in \Comp (\gG)$ satisfying the condition that $v(g,p,i) \in V(Y)$ for some $i=1,\dots ,\theta (p)$.
Then by Lemma \ref{lem:combine round-path} and Setting \ref{setting2} we see that the support $\mathrm{supp}(\eta_{\gG, p})$ is included in the $\delta_0$-neighborhood $B(\Cyl(p),\delta_0)$ of $\Cyl (p)$.

Note that by the definition of $\mathcal{O}$ we have
\begin{align*}
&\{ Y\in \Comp (\gG )\mid \ell (Y)\cap \bigcup_{h\in V(B_G(\id, r_0 ))}h\F \not=\emptyset \}\\
=&\bigsqcup_{p\in \mathcal{O}}\{ Y\in \Comp (\gG ) \mid v(g,p )\in V(Y) \text{ for }g\text{ with }p\in \R_r(g) \}.
\end{align*}
In other words, $v(g,p)\in V(Y)$ for $p\in \mathcal{O}\cap \R_r(g)$ if and only if
\[ \ell(Y)\cap \bigsqcup_{h\in V(B_G(\id, r_0 ))}h\Fla \not=\emptyset.\]
Since $K$ is included in $\bigsqcup_{h\in V(B_G(\id, r_0 ))}h\Fla$, we can assume that
if $\ell(Y)$ does not intersect $\bigsqcup_{h\in B_G(r_0,\id )}h\F$, then $\ell (Y)(\infty)$ does not belong to $A(K)$.
Hence, the following equality holds
\[ \eta_\gG|_{A(K)}=\sum_{p\in \mathcal{O}}\eta_{\gG,p}|_{A(K).}\]
As a result, by defining $\nu_p$ to be $\frac{1}{M} \eta_{\gG,p}$ for each $p\in \mathcal{O}$, we have
\[ \nu|_{A(K)}=\sum_{p\in \mathcal{O}}\nu_p|A(K).\]
This is the required equality in Condition (1) in Lemma \ref{lem:approximation data}.

(Condition (3) in Lemma \ref{lem:approximation data})
By the definition of $\eta_p$ for each $p\in \mathcal{O}$, we have
\begin{align*}
|\nu_p |
&=\frac{1}{M}\#\{Y\in \Comp (\gG )\mid v(g,p)\in V(Y)\}\\
&=\frac{1}{M}\#\{ v(g,p,i)\in V(\gG )\mid i=1,\dots ,\theta (p)\}\\
&=\frac{1}{M}\theta(p).
\end{align*}
Hence, for each $p\in \mathcal{O}$,
\[ \big| |\nu_p |-\mu (\Cyl (p))\big| =\left| \frac{1}{M}\theta(p)- \mu (\Cyl (p))  \right| <\hat{\varepsilon},\]
which is the required inequality.

Finally, since $\eta_{\gG,p}$ is a finite measure for every $p\in \mathcal{O}$, $\eta_\gG (A(K))$ is finite. We can assume that $K$ is sufficiently large.
Then we see that $\eta_\gG$ is a locally finite measure from Lemma \ref{lem:locally finite A(K) finite}.
Moreover, from Remark \ref{rem:discrete subset current}, $\eta_\gG$ is a discrete geodesic current.
Therefore, it follows by Lemma \ref{lem:approximation data} that $\nu=\frac{1}{M} \eta_\gG$ is a discrete geodesic current belonging to $U(\varepsilon; f_1,\dots ,f_\ell )$. \hspace{\fill}Q.E.D.

\section{Proof of denseness property of rational subset currents}\label{sec:dense subset currents}

Let $\gS$ be a cusped hyperbolic surface. Since our strategy for the proof of the denseness property of rational subset currents is the same as in the case of geodesic currents, we only present a sketch of the proof in this section.
We introduce the notion of an appropriate set of round-paths and the subset cylinder with respect to it, which plays the same role as a round-path and the cylinder with respect to it.

We use the setting in the beginning of Section \ref{sec:proof dense gc}, which we used in order to define the notion of round-paths.
Take the fundamental domain $\F$, the set $\mathrm{Side}(\F)$ of side-pairing transformations of $\F$ and the Cayley graph $\Cay(G)$ of $G$ with respect to the basis $\mathrm{Side}(\F)$ in the same manner as in Section \ref{sec:proof dense gc}.
We also take some horocycle parameter $\lambda$ and some radius $r\in \NN$.

\begin{definition}[Weak convex hull, set of round-paths and subset cylinder]
For $S\in \H (\partial \HH)$, we define the \ti{weak convex hull} $WCH(S)$ of $S$ to be the union of all geodesic lines connecting two points of $S$.

Let $g\in G$. Let $T$ be a set of round-paths of $B_G(g\F_\lambda, r)$, which includes some round-paths not passing through $g\Fla$ and is a finite set.
We say that $T$ is \ti{appropriate} if there exists $S\in \H (\partial \HH)$ satisfying the following two conditions:
\begin{enumerate}
\item for every round-path $p\in T$, there exists a geodesic line $\ell$ connecting two points of $S$ such that $\ell $ passes through $p$, i.e., for $p=[e_1,\dots ,e_k]\in T$, $\ell$ passes through $e_1,\dots ,e_k$ in this order while passing through $\overline{B_g(g\F_\lambda,r)}$;
\item for every geodesic line $\ell$ connecting two points of $S$, there exists $p\in T$ such that $\ell$ passes through $p$.
\end{enumerate}
If $S\in \H (\partial \HH )$ satisfies the above two conditions, we say that the restriction of $WCH(S)$ to $B_G(g\F_\lambda, r)$ equals $T$.

For an appropriate set $T$ of round-paths, we define the \ti{subset cylinder} $\SCyl (T)$ with respect to $T$ to be the subset of $\H (\partial \HH )$ consisting of an element $S$ satisfying the condition that the restriction of $WCH(S)$ to $B_G(g\F_\lambda, r)$ equals $T$.
We denote by $\R^\ast_r(g)$ the set of all appropriate sets of round-paths that contains a round-path containing an edge $e$ of $g\F_\lambda$ with $e\cap g\F\not=\emptyset$. Note that $\R^\ast_r(g)$ is a finite set.
\end{definition}

The notion of subset cylinder has the same property as that of cylinder and we have the following equality
\[ \bigsqcup_{T\in \R^\ast_r(g)}\SCyl (T)=\{ S\in \H (\partial \HH )\mid CH(S)\cap g\F_\lambda\not=\emptyset \} =A(g\F_\lambda ).\]
The properties that we proved in Section \ref{sec:proof dense gc} from Lemma \ref{lemma:diam of cyl} to Lemma \ref{lem:approximation data} except Lemma \ref{lemma:diam quasi-geod} can be naturally generalized to the subset current version.

The subset current version of Lemma \ref{lemma:diam quasi-geod} is as follows:
\begin{lemma}\label{lem:diam weak convex hull}
Fix $g\in G$, $a\geq 1, b\geq 0$ and $\delta_0>0$.
There exist a large horocycle parameter $\lambda$ and a constant $\delta_1>0$ such that for a sufficiently large $r\in \NN$ and any $T\in \R_r^\ast(g)$, if a set $Y$ of $(a,b)$-quasi-geodesics satisfies the following conditions:
\begin{enumerate}
\item $|Y|=\bigcup_{\ell\in Y} \ell$ is $\delta_1$-quasi-convex;
\item for every $\ell \in Y$ passing through $\overline{B_G(g\F_\lambda ,r)}$, there exists $p\in T$ such that $\ell$ passes through $p$;
\item for every $p\in T$, there exists $\ell \in Y$ such that $\ell$ passes through $p$,
\end{enumerate}
then the limit set $|Y|(\infty )$ of $|Y|$ is contained in the $\delta_0$-neighborhood of $\SCyl(T)$.
\end{lemma}
\begin{proof}
Recall that $|Y|$ is $\delta_1$-quasi-convex if for any $x,y\in |Y|$, the geodesic segment $[x,y]$ connecting $x$ to $y$ is included in the $\delta_1$-neighborhood of $|Y|$ with respect to the hyperbolic metric on $\HH$.

To obtain a contradiction, we suppose that $\gL (|Y|)$ is not contained in the $\delta_0$-neighborhood of $\SCyl(T)$.
Then there exists $\delta_0'>0$ depending on $\delta_0$ such that we can take $S\in \SCyl (T)$ and $x \in S$ such that
\[ B_{\ol{\HH}}(x ,\delta_0') \cap WCH(|Y|(\infty ))=\emptyset ,\]
or we can take $x\in |Y|(\infty )$ and $S\in \SCyl (T)$ such that
\[ B_{\ol{\HH}}(x ,\delta_0') \cap WCH(S)=\emptyset.\]

First, we consider the former case. By the same argument as in Lemma \ref{lemma:diam of cyl}, if $\lambda, r$ are sufficiently large, then there exists a round-path $p=[e_1,\dots ,e_k]$ such that $e_k$ is sufficiently close to $x$ in $\ol{\HH}$, which implies that one of the end-points of an $(a,b)$-quasi-geodesic line $\ell$ passing through $p$ must contained in $B_{\ol{\HH}}(x,\delta_0'/2 )$. Therefore, we can see that
\[ B_{\ol{\HH}}(x ,\delta_0') \cap WCH(|Y|(\infty ))\not=\emptyset ,\]
which is a contradiction.

Next, we consider the latter case. In this case, the quasi-convexity of $|Y|$ plays an important role because $x$ can be an end-point of a quasi-geodesic far away from $B_G(g\F_\lambda ,r)$.
Take $\ell \in Y$ such that one of the end-points of $\ell$ is sufficiently close to $x$ and take $\ell_1\in Y$ passing through some round-path $p_1\in T$ containing an edge of $g\F_\lambda$.
Since $|Y|$ is $\delta_1$-quasi-convex, by considering a geodesic segment connecting $y\in \ell \cap B_{\ol{\HH}}(x ,\delta_0'/4)$ to some point on $\ell_1$, which is included in the $\delta_1$-neighborhood of $|Y|$, there exists $\ell_2\in Y$ passing through $p_2\in T$ such that $p_2$ contains an edge $e$ of $B_G(g\F_\lambda, r)$ sufficiently close to $x$ in $\ol{\HH}$.
Then $WCH(S)$ includes a geodesic line passing through $p_2$, which must
intersect $B_{\ol{\HH}}(x ,\delta_0')$, i.e.,
\[B_{\ol{\HH}}(x ,\delta_0') \cap WCH(S)\not=\emptyset,\]
which is a contradiction.
\end{proof}

The notion of connectability can be also generalized to the subset current version.
For adjacent $u,v\in V(\Cay (G))$, $T_1\in \R_r^\ast (u)$ and $T_2\in \R_r^\ast (v)$, we say that $T_1$ and $T_2$ are \ti{connectable} if the following conditions follows:
\begin{enumerate}
\item $T_1$ contains a round-path containing an edge $e$ of $v\Fla$ with $e\cap v\Fla \neq \emptyset$;
\item $T_2$ contains a round-path containing an edge $e$ of $u\Fla$ with $e\cap u\Fla \neq \emptyset$;
\item the restriction of $T_1$ to $B_G(u\F_\lambda ,v\F_\lambda, r)$ coincides with the restriction of $T_2$ to $B_G(u\F_\lambda ,v\F_\lambda, r)$.
\end{enumerate}

Fix non-zero measure $\mu \in \SC (\gS)$. Then we can approximate $\mu$ by
\[ \theta\: \bigsqcup_{g\in G} \R_r^\ast (g)\rightarrow \mathbb{Z}_{\geq 0}\]
satisfying the subset current version of the conditions in Section \ref{sec:proof dense gc}.
From $\theta$, we construct a graph $\gG$ that $G$ acts on in the same manner as in Section \ref{sec:proof dense gc}, i.e.,
\[ V(\gG )=\{ v(g,T,i )\}_{g\in G, T\in \R_r^\ast (g), i=1,\dots , \theta (T)} \]
and if $v(u,T)$ and $v(v,T')$ is connected by an edge, then $u,v$ are adjacent in $\Cay (G)$, and $T$ and $T'$ are connectable.

Now, we consider each connected component $Y$ of $\gG$ and construct a set of quasi-geodesics by combining round-paths
\[ \{ p\}_{v(g,T)\in V(Y), p\in T}.\] The biggest difference between the case of subset currents and the case of geodesic currents is that $Y$ is a sub-tree of $\Cay (G)$, which is much more complicated than a finite segment or a bi-infinite line.
Define $\Comp (\gG)$ to be the set of all connected components of $\gG$.

Let $Y\in \Comp (\gG)$, $v(g,T)\in V(Y)$ and $p\in T$.
We remark that $p$ may not passes through an edge of $g\F_\lambda$ but passes through $h_0\F_\lambda$ for some $h_0\in V(B_G(g,r))$.
Then we can take a geodesic path of vertices in $Y$ connecting $v(g,T)$ to $v(h_0,T_0)$ since there exists $S\in \H (\partial \HH)$ such that the restriction of $WCH(S)$ to $B_G(g\F_\lambda ,r)$ equals $T$, which must passes through every $g'\F_\lambda$ for every $g'\in G$ on a geodesic path connecting $g$ to $h_0$ in $\Cay (G)$. In addition, $T_0$ contains a round-path $p_0$ including $p$ as a sub-round-path.

Considering the extension of $p_0$ by the connectability, we can obtain a finite or bi-infinite sequence $\{ v(h_i,T_i)\}_i$ of $V(Y)$ centered at $v(h_0,T_0)$ satisfying the following conditions
\begin{enumerate}
\item $v(h_i,T_i)$ is connected to $v(h_{i+1},T_{i+1})$ by an edge for every $i$ (except the case in which the vertex $v(h_{i+1},T_{i+1})$ does not exist);
\item for each $i$, there exists $p_i\in T_i$ such that $p_i$ and $p_{i+1}$ are connectable;
\item the sequence $\{ v(h_i,T_i)\}_i$ ends at $i_1$ or $-i_2$ for $i_1,i_2\in \mathbb{Z}_{\geq 0}$ if and only if one of the end edges of $p_{i_1}$ or $p_{-i_2}$ is a horocyclic edge.
\end{enumerate}
Note that even when $\{ v(h_i,T_i)\}$ is infinite, we do not know whether it is bi-infinite or not.

For the sequence $\{ v(h_i,T_i)\}$, by combining the round-paths $\{p_i\}$, we can take a sequence $[e_i]$ of edges, which is independent of the choice of $h_0$.
Then in the same manner as in Lemma \ref{lem:combine round-path} we can obtain an $(a,b)$-quasi-geodesic $\ell(Y,p)$ passing through $[e_i]$ in this order. Once we fix $\ell(Y,p)$, we define $\ell(Y,p')$ to be $\ell(Y,p)$ for every $p'\in T'$ for $v(u,T')\in V(Y)$ satisfying the condition that $p'$ is a sub-sequence of the sequence $[e_i]$.
We define $|Y|$ as
\[ |Y|:=\bigcup_{v(g,T)\in V(Y), p\in T}\ell(Y,p). \]

By the definition, we see that $\{\ell (Y,p) \}$ satisfies the condition 2 and 3 in Lemma \ref{lem:diam weak convex hull}.
We check that $|Y|$ is $\delta_1$-quasi-convex for some constant $\delta_1$, which depends on $a, b,\lambda, \HH$ and a constant $s$ in the following lemma.

\begin{lemma}\label{lem:quasi convex and horodisk}
Let $v(g_{-m},T_{-m}),\dots, v(g_m,T_m)$ be a geodesic path of vertices in $Y$ for $m\in \NN$.
Let $H$ be a connected component of $\HH\setminus \HH_\lambda$, which is a horodisk centered at a parabolic fixed point $\xi$ of $G$.
Assume that for each $i\in \{ -m ,\dots, m\}$ a horocyclic edge $e_i$ of $g_i\F_\lambda$ is included in the boundary $\partial H$ of $H$.
Then there exists $s\in \NN$ depending only on $\lambda$ such that if $m\geq 3s$, then $|Y|$ includes a geodesic ray emanating from one of $e_{-m},\dots ,e_m$ to $\xi$.
We also assume that the radius $r$ for sets of round-paths is much larger than $s$.
\end{lemma}
\begin{proof}
Consider the upper-half plane model of $\HH$ and assume that $H$ is a horodisk centered at $\infty$.
Then the boundary $\partial H$ of $H$ is $\{ x+y_0\sqrt{-1} \mid x\in \RR\}$ for some $y_o\in \RR$ depends on $\lambda$.
The endpoints of $e_0$ are $x_0+y_0\sqrt{-1}$ and $x_0+\alpha +y_0\sqrt{-1}$ for some $\alpha>0$.
We can assume that $x_0+y_0\sqrt{-1}$ is the right endpoint of $e_{-1}$ and $x_0+\alpha +y_0\sqrt{-1}$ is the left endpoint of $e_1$.
Let $s$ be the smallest positive integer satisfying the condition that $s\alpha>y_0+2\alpha$.

To obtain a contradiction, suppose that $m\geq 3s$ and $|Y|$ does not include any geodesic rays emanating from one of $e_{-m},\dots ,e_m$ to $\infty$.
Take $S\in \SCyl (T_0)$.
Since $r>s$, there exist $\xi_i,\zeta_i \in S$ such that the geodesic $[\xi_i ,\zeta_i]$ passes through $g_i\Fla$ for $i=-s,\dots, s$. By the assumption, $\xi_i$ and $\zeta_i$ belong to $\RR$ (we assume that $\xi_i<\zeta_i$), and $[\xi_i,\zeta_i]$ does not intersect two non-adjacent edges of $e_{-s},\dots, e_s$, which implies that $\zeta_i-\xi_i< 2y_0+2\alpha$ for $i=-s, \dots ,s$ (see Figure \ref{fig:quasi-convex-lemma}).
\begin{figure}[h]
\centering 
\includegraphics{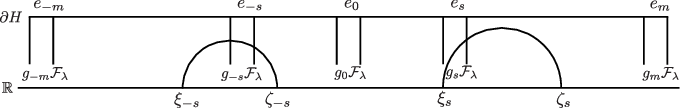}
\caption{Setting of the proof of Lemma \ref{lem:quasi convex and horodisk}.}\label{fig:quasi-convex-lemma}
\end{figure}

Then we can see that
\[ \xi_i<x_0+(i+1)\alpha \text{ and } \zeta_i >x_0+i\alpha\]
since $[\xi_i,\zeta_i]\cap g_i\Fla\neq \emptyset$.
Hence,
\[ \zeta_s-\xi_{-s}>x_0+s\alpha -(x_0+(-s+1)\alpha )=2s\alpha -\alpha >2y_0+2\alpha. \]
Note that
\[ \zeta_s<\xi_s+2y_0+2\alpha <x_0+(s+1)\alpha +2s\alpha < x_0+(3s+1)\alpha \leq x_0+(m+1)\alpha \]
and
\[ \xi_{-s}>\zeta_{-s}-2y_0-2\alpha >x_0-s\alpha-2s\alpha >x_0+(-3s)\alpha \geq x_0+(-m)\alpha.\]
Therefore, $[\xi_{-s},\zeta_{s}]$ intersects with two non-adjacent edges of $e_{-m},\dots, e_m$. This implies that $|Y|$ includes a geodesic ray emanating from one of $e_{-m},\dots ,e_m$ to $\infty$, which is a contradiction.
\end{proof}
\begin{proof}[Proof of the quasi-convexity of $|Y|$]
Let $x,y\in |Y|$. Take $v(g,T),v(g',T')\in V(Y)$ and $p\in T, p'\in T'$ such that $x\in \ell (Y,p)$ and $y\in \ell (Y,p')$.
When $x$ (or $y$) belongs to a geodesic ray in a horodisk of $\HH\setminus \HH_\lambda$, it is enough to consider the nearest point $x'\in \HH_\lambda\cap \ell (Y,p)$ from $x$ (or the nearest point $y'\in \HH_\lambda \cap \ell (Y,p')$ from $y$).
Hence, we can assume that $x\in g\Fla$ and $y\in g'\Fla$.

Take a geodesic path of vertices
\[v(g_0,T_0)=v(g,T), v(g_1,T_1),\dots ,v(g_k,T_k)=v(g',T')\]
in $Y$.
From the shape of the fundamental domain $\F$, we can see that the geodesic $[x,y]$ passes through $g_0\F, g_1\F, \dots ,g_k\F$ in this order.
If $[x,y]$ is included in $\HH_\lambda$, then $[x,y]$ is included in the $\mathrm{diam}(\Fla)$-neighborhood of $|Y|$ since for every $g_i\Fla$, there exists a quasi-geodesic of $|Y|$ passes through $g_i\Fla$.

Now, we consider the case in which $[x,y]$ is not included in $\HH_\lambda$.
Assume that $[x,y]$ goes into a horodisk $H$ while passing through $g_i\F$ and goes out from $H$ while passing through $g_{i+t}\F$.
If $t<6s$, then there exists a constant $C_s>0$ depending on $a,b,\lambda ,s$ such that
\[ [x,y]\cap \bigcup_{j=i}^{i+t}g_j\F\]
is included in the $(\mathrm{diam}(\Fla)+C_s)$-neighborhood of $|Y|$.
If $t\geq 6s$, then $[x,y]$ intersects some geodesic rays that are pieces of quasi-geodesics of $|Y|$ while passing through $\bigcup_{j=i}^{i+t}g_j\F$ by Lemma \ref{lem:quasi convex and horodisk}.
Therefore, $[x,y]$ is included in the $\delta_1$-neighborhood of $|Y|$ for a constant $\delta_1$ depending on $a,b,\lambda$ and $s$.
\end{proof}

Now, we can obtain a discrete subset current
\[ \eta_\Gamma :=\sum_{Y\in \Comp (\gG)}\delta_{|Y|(\infty )} \in \SC (\gS ),\]
and we can prove that $\frac{1}{M} \eta_\gG$ approximates the given $\mu \in \SC (\gS)$ for some $M\in \NN$.
We omit the rest of the proof since it is almost the same as in the case of geodesic currents.

\end{document}